\documentclass[a4paper]{amsart}

\usepackage{amssymb}
\usepackage[all]{xypic}

\def\frak{\mathfrak}
\def\Bbb{\mathbb}
\def\Cal{\mathcal}

\swapnumbers
\numberwithin{equation}{subsection}

\newtheorem{prop}[subsection]{Proposition}
\newtheorem*{prop*}{Proposition}
\newtheorem{thm}[subsection]{Theorem}
\newtheorem*{thm*}{Theorem}
\newtheorem{lem}[subsection]{Lemma}
\newtheorem*{lem*}{Lemma}

\newtheorem*{kor*}{Corollary}

\newcommand{\ad}{\operatorname{ad}}
\newcommand{\Ad}{\operatorname{Ad}}
\renewcommand{\exp}{\operatorname{exp}}

\newcommand{\Fl}{\operatorname{Fl}}

\newcommand{\x}{\times}
\renewcommand{\o}{\circ}

\let\ccdot\cdot
\def\cdot{\hbox to 2.5pt{\hss$\ccdot$\hss}}

\newcommand{\al}{\alpha}

\newcommand{\ka}{\kappa}
\newcommand{\la}{\lambda}
\newcommand{\om}{\omega}
\renewcommand{\phi}{\varphi}
\newcommand{\ph}{\varphi}
\newcommand{\si}{\sigma}

\newcommand{\ze}{\zeta}

\newcommand{\Ga}{\Gamma}
\newcommand{\La}{\Lambda}
\newcommand{\Ph}{\Phi}
\newcommand{\Ps}{\Psi}
\newcommand{\Rho}{{\mbox{\sf P}}}
\newcommand{\Om}{\Omega}

\renewcommand{\Up}{\Upsilon}
\newcommand{\cg}{{\Cal G}}
\newcommand{\ca}{{\Cal A}}

\begin{document}
\title{Weyl structures for parabolic geometries} 
\author{Andreas \v Cap and Jan Slov\'ak}

\address{A.C.: Institut f\"ur Mathematik, Universit\"at Wien,
Strudlhofgasse~4, A--1090 Wien, Austria, and International Erwin
Schr\"odinger Institute for Mathematical Physics, Boltzmanngasse 9, A-1090
Wien, Austria\newline\indent J.S.: Department of Algebra and Geometry,
Masaryk University, Jan\'a\v ckovo n\'am. 2a, 662~95~Brno, Czech Republic}

\email{Andreas.Cap@esi.ac.at, slovak@math.muni.cz}

\begin{abstract}
Motivated by the rich geometry of conformal Riemannian manifolds and
by the recent development of  geometries modeled on homogeneous spaces
$G/P$ with $G$ semisimple and $P$ parabolic,  Weyl structures and
preferred 
connections are introduced in this general framework. In particular, we
extend the notions of scales, closed and exact Weyl connections, and
Rho--tensors, we characterize the classes of such objects, and we use the
results to give a new description of the Cartan bundles and connections for all
parabolic geometries.
\end{abstract}

\maketitle
\section{Introduction}\label{1}

Cartan's generalized spaces are curved analogs of the homogeneous
spaces $G/P$ defined by means of an absolute parallelism on a principal
$P$--bundle. This very general framework was originally built in connection
with the equivalence problem and Cartan's general method for its
solution, cf. e.g. \cite{Cartan-conf}. Later on, however, these ideas got
much more attention. In particular, several well known geometries were
shown to allow a canonical object of that type with suitable choice of
semisimple $G$ and parabolic $P$, see e.g. \cite{Kob}. Cartan's original 
approach was generalized and extended for all such groups, cf. 
\cite{Tan, Mor, Yam, CS}, and links to other areas were discovered, see e.g.
\cite{BE, Bas, CheM}. The best known examples are the conformal Riemannian, 
projective, almost quaternionic, and CR structures and the common name
adopted is {\em parabolic geometries\/}.

The relation to twistor theory renewed the interest in a good
calculus for such geometries, which had to improve the techniques in
conformal geometry and to extend them to other geometries. 
Many steps in this direction were done, see for 
example \cite{Tho1,Tho2, Wuensch, Gov2} for classical methods 
in conformal geometry, and \cite{BEG, BaiE, Gov1, Gov3, Gov4} for
generalizations. 

A new approach to this topic, motivated mainly by \cite{Och, Bas, BE}, was
started in \cite{CSS1-3, CSS4}. The novelty consists in the combination of
Lie algebraic tools with the frame bundle approach to all objects and we
continue in this spirit here. Our general setting for Weyl structures and 
scales has been also inspired by \cite{BaiE, Gau}.

In Section \ref{2} we first outline some general aspects of parabolic
geometries and then we present the basic objects like tangent and cotangent
bundles and the curvature of the geometry in a somewhat new perspective.
This will pave our way to the Weyl structures in the rest of the paper. Our
basic references for Section \ref{2} are \cite{CS} and \cite{Sl-notes}, the
reader may also consult \cite{CSS4}.  For the classical point of view of
over--determined systems, we refer to \cite{Tan, Yam} and the references 
therein.

The Weyl structures are introduced in the beginning of Section \ref{3}.
Exactly as in the conformal Riemannian case, the class of Weyl structures
underlying a parabolic geometry on a manifold $M$ is always an affine space
modeled on one--forms on $M$ and each of them determines a linear
connection on $M$. Moreover, the difference between the linear connection
induced by a Weyl structure and the canonical Cartan connection is encoded
in the so called Rho--tensor (used heavily in conformal geometry since
the beginning of the century).  Next, we define the bundles of scales as
certain affine line bundles generalizing the distinguished bundles of
conformal metrics, and we describe the correspondence between connections on
these line bundles and the Weyl structures, see Theorem \ref{3.12}.  On the
way, we achieve explicit formulae for the deformation of Weyl structures and
the related objects in Proposition \ref{3.9}, which offers a generalization
for the basic ingredients of various calculi. The exact Weyl geometries are
given by scales, i.e. by (global) sections of the bundles of scales, thus
generalizing the class of Levi--Civita connections for conformal geometries.
At the same time, this point of view leads to a new presentation of the
canonical Cartan bundle as the bundle of connections on the bundle of scales
(pulled back to the defining infinitesimal flag structure, cf. \ref{2.7} and
\ref{3.12}). 
In the end of Section \ref{3}, we
define another class of distinguished local Weyl structures which achieve
the best possible approximation of the canonical Cartan connections, see
Theorem \ref{3.16}. In the
conformal case, these normal Weyl structures improve the construction of the
Graham's normal coordinates, cf. \cite{Gra}.

The last section is devoted to characterizations of all the
objects related to a choice of a Weyl structure. More explicitly, the
ultimate goal is to give a recipe how to decide which soldering forms and
linear connections on a manifold $M$ equipped with a regular infinitesimal
flag structure are obtained from a Weyl--structure and to compute the
corresponding Rho--tensor. For
this purpose, we define the general Weyl forms and their Weyl curvatures and
the main step towards our aim is achieved in Theorem \ref{4.4}. Next, we
introduce the total curvature of a Weyl form which is easier to interpret on
the underlying manifold than the Weyl curvature. The characterization is
then obtained by carefully analyzing the relation between these two
curvatures.

This entire paper focuses on the introduction of new structures and their
nice properties. We should like to mention that essential use of these new
concepts has appeared already in \cite{CSS4} and \cite{Cal}. 

\medskip 
\noindent{\bf Acknowledgements.} 
The initial ideas for this research evolved during the stay of the second
author at the University of Adelaide in 1997, supported by the Australian
Research Council. The final work and writing was done at the Erwin
Schr\"odinger Institute for Mathematical Physics in Vienna. The second
author also acknowledges the support from GACR, Grant Nr. 201/99/0296.
Our thanks are also due to our colleagues for many discussions.

\section{Some background on parabolic geometries}\label{2}

\subsection{$|k|$--graded Lie algebras}\label{2.1}
Let $G$ be a real or complex semisimple Lie group, whose Lie algebra
$\frak g$ is equipped with a grading
of the form 
$$
\frak g=\frak g_{-k}\oplus\dots\oplus\frak g_0\oplus\dots\oplus\frak g_{k}
.$$
Such algebras $\frak g$ are called {\em $|k|$--graded Lie algebras}. 

Throughout this paper we shall further assume that  no simple ideal of $\frak
g$ is contained in $\frak g_0$ and that the (nilpotent) subalgebra $\frak
g_-=\frak g_{-k}\oplus\dots\oplus\frak g_{-1}$ is generated by $\frak
g_{-1}$. Such algebras are sometimes called {\em effective semisimple graded
Lie algebras of $k$-th type}, cf. \cite{Kan, Tan}. By $\frak p_+$ we denote
the subalgebra $\frak g_1\oplus\dots\oplus\frak g_k$ and by $\frak p$ the
subalgebra $\frak g_0\oplus\frak p_+$. We also write $\frak g_-=\frak
g_{-k}\oplus\dots\oplus\frak g_{-1}$, and $\frak g^j=\frak g_j\oplus\dots\oplus
\frak g_k$, $j=-k,\dots,k$.

It is well known that then $\frak p$ is a parabolic subalgebra of
$\frak g$, and actually the grading is completely determined by this
subalgebra, see e.g.~\cite{Yam}, Section 3. Thus all complex simple
$|k|$--graded $\frak g$ are classified by subsets of simple roots of complex
simple Lie algebras (i.e. arbitrary placement of crosses over the Dynkin
diagrams in the notation of \cite{BE}), up to conjugation. 
The real $|k|$--graded simple Lie algebras are classified easily by means of 
Satake diagrams: the $|k|$--grading of the complex simple $\frak g$ induces a
$|k|$--grading on a real form if and only if (i) only `white' nodes in the
Satake diagram have been crossed out, and, (ii) if a node is crossed out,
then all nodes connected to this one by the double arrows in the Satake
diagram have to be crossed out too, see \cite{Kan} or \cite{Yam} for  more 
details. Very helpful notational conventions
and computational recipes can be found in \cite{BE}. 

\subsection{}\label{2.2}
Let us recall basic properties of Lie groups $G$ with 
(effective) $|k|$--graded Lie algebras $\frak g$.
First of all, there is always a unique element $E\in \frak g_0$ with the
property $[E,Y]=jY$ for all $Y\in \frak g_j$, $j=-k,\dots, k$, the {\em
grading element}. Of course, $E$ belongs to the center $\frak z$ of the
reductive part $\frak g_0$ of $\frak p\subset \frak g$. 

The Killing form provides isomorphisms $\frak g_i^*\simeq \frak g_{-i}$
for all $i=-k,\dots,k$ and, in particular, its restrictions to the center
$\frak z$ and the semisimple part $\frak g_0^{ss}$ of $\frak g_0$ are
non--degenerate.

Now, there is the closed subgroup $P\subset G$ of all elements whose adjoint
actions leave the $\frak p$--submodules $\frak g^j=\frak
g_j\oplus\dots\oplus \frak g_k$ invariant, $j=-k,\dots,k$. The Lie algebra
of $P$ is just $\frak p$ and there is the subgroup $G_0\subset P$ of
elements whose adjoint action leaves invariant the grading by $\frak
g_0$--modules $\frak g_i$, $i=-k,\dots,k$. This is the reductive part of the
parabolic Lie subgroup $P$, with Lie algebra $\frak g_0$. 
We also define subgroups $P_+^j=\operatorname{exp}(\frak
g_j\oplus\dots\oplus \frak
g_k)$, $j=1,\dots,k$, and we write $P_+$ instead of $P_+^1$. Obviously
$P/P_+=G_0$  and $P_+$ is nilpotent. Thus $P$ is the
semidirect product of $G_0$ and the nilpotent part $P_+$. More
explicitly, we have (cf. \cite{CS}, Proposition 2.10, or \cite{Tan, Yam})

\begin{prop}\label{2.3} 
For each element $g\in P$, there exist unique elements $g_0\in G_0$ and $Z_i\in
\frak g_i$, $i=1,\dots,k$, such that
$$
g=g_0\operatorname{exp}Z_1
\operatorname{exp}Z_2\dots\operatorname{exp}Z_k
.$$
\end{prop}

\subsection{Parabolic geometries}\label{2.4}
Following Elie Cartan's idea of generalized spaces (see \cite{Sha} for a
recent reading), a curved analog of the homogeneous space $G/P$ is a
right invariant absolute parallelism $\om$ on a principal $P$--bundle $\Cal
G$ which reproduces the fundamental vector fields. In our approach, a 
(real) {\em parabolic geometry $(\Cal G,\om)$
of type $G/P$} is a principal fiber bundle $\Cal G$ with structure group
$P$, equipped with a smooth one--form $\om\in\Om^1(\Cal G,\frak g)$ satisfying
\begin{enumerate}
\item[(1)] $\om(\ze_Z)(u)=Z$ for all $u\in\Cal G$ and fundamental fields $\ze_Z$,
$Z\in \frak p$
\item[(2)] $(r^b)^*\om = \Ad(b^{-1})\o \om$ for all $b\in P$
\item[(3)] $\om|_{T_u\Cal G}: T_u\Cal G\to \frak g$ is a linear isomorphism for
all $u\in \Cal G$.
\end{enumerate} 
In particular, each $X\in \frak g$ defines the {\em constant vector field}
$\om^{-1}(X)$ defined by $\om(\om^{-1}(X)(u))=X$, $u\in \Cal G$. In this
paper, we shall deal with smooth real parabolic geometries only. The one
forms with properties (1)--(3) are called {\em Cartan connections}, cf.
\cite{Sha}.

The morphisms between parabolic geometries $(\cg,\om)$ and $(\cg',\om')$ are
principal fiber bundle morphisms $\ph$ which
preserve the Cartan connections, i.e. $\ph:\cg\to \cg'$ and $\ph^*\om'=\om$. 

\subsection{The curvature}\label{2.5} 
The structure equations define the horizontal smooth form 
$K\in\Om^2(\Cal G, \frak g)$ called the {\em curvature} of the Cartan
connection $\om$:
$$
d\om + \frac 12[\om,\om] = K
.$$
The {\em curvature function} $\ka:\Cal G\to \wedge^2\frak g_{-}^*\otimes
\frak g$ is then defined by means of the
parallelism 
$$
\ka(u)(X,Y)=K(\om^{-1}(X)(u),\om^{-1}(Y)(u))= 
[X,Y] - \om([\om^{-1}(X),\om^{-1}(Y)])
.$$
In particular, the curvature function is valued in the cochains for the
second cohomology $H^2(\frak g_{-},\frak g)$. Moreover, there are two ways
how to split $\ka$. We may consider the target components $\ka_i$
according to the values in $\frak g_i$. The whole $\frak g_-$--component 
$\ka_-$ is called the {\em torsion} of the Cartan connection $\om$.
The other possibility is to consider the homogeneity of the bilinear maps
$\ka(u)$, i.e. 
$$
\ka=\sum_{\ell=-k+2}^{3k}\ka^{(\ell)},\quad \ka^{(\ell)}:\frak g_i\times\frak
g_j\to \frak g_{i+j+\ell} 
.$$

Since we deal with semisimple algebras only, there is the
codifferential $\partial^*$ which is ajoint to the Lie algebra cohomology 
differential
$\partial$, see e.g. \cite{Kostant}. Consequently, there is the Hodge theory
on the cochains which enables to deal very effectively with the curvatures.
In particular, we may use several restrictions on the values of the
curvature which turn out to be quite useful.

\subsection{Definition}\label{2.6}
The parabolic geometry $(\Cal G,\om)$ with the curvature function $\ka$ is
called {\em flat} if $\ka=0$, {\em torsion--free} if $\ka_-=0$, 
{\em normal} if $\partial^*\o\ka=0$, and {\em regular} if it is normal and
$\ka^{(j)}=0$ for all $j\le 0$. 

\smallskip

Obviously, the morphisms of parabolic geometries preserve the above types
and so we obtain the corresponding full subcategories of regular,  normal,
torsion free, and flat parabolic geometries of a fixed type $G/P$. See
\cite{CSS4}, Section 2, for more details.

\subsection{Flag structures}\label{2.7}
The homogeneous models for parabolic geometries are the real generalized
flag manifolds $G/P$. Curved parabolic geometries look like $G/P$
infinitesimally. Indeed, the filtration of
$\frak g$ by the $\frak p$--submodules $\frak g^j$ is transfered to the
right invariant filtration $T^{j}\Cal G$ on the tangent space $T\Cal G$ by
the parallelism $\om$. The tangent projection $Tp:T\Cal G\to TM$ 
then provides the filtration $TM=T^{-k}M\supset T^{-k+1}M\supset\dots\supset
T^{-1}M$ of the tangent space of the underlying manifold $M$. Moreover,
the structure group of the associated graded tangent space
$\operatorname{Gr}TM=(T^{-k}M/T^{-k+1}M)\oplus \dots \oplus
(T^{-2}M/T^{-1}M)\oplus T^{-1}M$ reduces  automatically to $G_0$ since $\Cal
G_0=\Cal G/P_+$ clearly plays the role of its frame bundle. 
The following lemma is not difficult to prove, 
see e.g. \cite{SchmSl}, Lemma 2.11.

\begin{lem*} Let $(\Cal G,\om)$ be a parabolic geometry, $\ka$ its curvature
function. Then $\ka^{(j)}=0$ for all $j<0$ if and only if the Lie bracket of
vector fields on $M$ is compatible with the filtration, i.e. $[\xi,\eta]$ is
a section of $T^{i+j}M$ for all sections $\xi$ of $T^iM$, and  $\eta$ of 
$T^{j}M$. Hence it
defines an algebraic bracket $\{\ ,\ \}_{\text{Lie}}$ 
on $\operatorname{Gr}TM$. Moreover, this bracket coincides with the
algebraic bracket $\{\ ,\ \}_{\frak g_0}$ defined on $\operatorname{Gr}TM$
by means of the $G_0$--structure if and only if $\ka^{(j)}=0$ for all $j\le
0$.
\end{lem*}

We call the filtrations of $TM$ with reduction of $\operatorname{Gr}TM$ to
$G_0$  satisfying the very last condition
of the lemma the {\em regular infinitesimal flag structures of type $\frak
g/\frak p$}. In fact, the structures clearly depend on the choice of the Lie
group $G$ with the given Lie algebra $\frak g$. This choice is always
encoded already in $G_0$. On the other hand, there are always several
distinguished choices, e.g. the full automorphism group of $\frak g$, the
adjoint group, and the unique connected and simply connected group. In the
conformal geometries these choices lead to conformal Riemannian manifolds,
oriented conformal menifolds, and (oriented) conformal spin manifolds,
respectively. Obviously, the various choices of $G$ do not matter much
locally and we shall not discuss them explicitly in this paper.

The $G_0$ structures on $\operatorname{Gr}TM$ are equivalent to the frame
forms of length one defined and used in \cite{CS} while the condition 
$\ka^{(j)}=0$ for all $j\le 0$ is equivalent to the structure equations 
for these frame forms imposed in the construction of \cite{CS}. In view of
this relation, we also call our bundles $\cg_0$ equipped with the regular
infinitesimal flag structures the {\em $P$--frame bundles of degree one}. 
In particular, we obtain (see \cite{CS}, Section 3)

\begin{thm}\label{2.8}
There is the bijective correspondence between the isomorphism classes of
regular parabolic geometries of type $G/P$ and the regular infinitesimal 
flag structures of type $\frak g/\frak p$ on $M$, except for one series of
one--graded, and one series of two--graded Lie algebras $\frak g$ for which 
$H^1(\frak g_-,\frak g)$ is nonzero in homogeneous degree one.
\end{thm} 

Both types of the exceptional geometries from the Theorem will be mentioned
in the examples below.

\subsection{Example}\label{2.9}
The parabolic geometries with $|1|$--graded Lie algebras $\frak g$ are
called {\em irreducible}. Their tangent bundles do not carry any nontrivial
natural filtration and this irreducibility of $TM$ is reflected in the name.
The classification of all such simple real Lie algebras is well known (cf.
\cite{KoNa} or \ref{2.1} above). We may list all the corresponding
geometries, up to the possible choices of the groups $G_0$, roughly as follows:
\begin{itemize}
\item[$A_{\ell}$] the split form, $\ell>2$ --- 
the {\em almost Grassmannian structures} with
homogeneous models of $p$--planes in $\Bbb R^{\ell+1}$, $p=1,\dots,\ell$.
The choice $p=1$ yields the projective structures which represent one of the
two exceptions in \ref{2.8}.
\item[$A_{\ell}$] the quaternionic form, $\ell=2p+1>2$ --- the {\em 
almost quaternionic
geometries} in dimensions $4p$, and more general geometries modeled on
quaternionic Grassmannians.
\item[$A_{\ell}$] one type of geometry for the algebra $\frak s\frak
u(p,p)$, $\ell = 2p-1$.
\item[$B_{\ell}$] the {\em (pseudo) conformal geometries} in all odd dimensions
$2m+1\ge 3$.
\item[$C_{\ell}$] the split form, $\ell > 2$ --- the {\em almost Lagrangian
geometries} modeled on the Grassmann 
manifold  of maximal Lagrangian subspaces in the symplectic $\Bbb R^{2\ell}$.
\item[$C_\ell$] another type of geometry corresponding to the algebra $\frak
s\frak p(p,p)$, $\ell=2p$.
\item[$D_\ell$] the {\em (pseudo) conformal geometries} in all even dimensions
$m\ge 4$.
\item[$D_\ell$] the real almost spinorial geometries with $\frak g=\frak s\frak
o(p,2\ell-p)$, $p=1,\dots, \ell-2$.
\item[$D_\ell$] the quaternionic almost spinorial geometries with $\frak
g=\frak u^*(\ell,\Bbb H)$.
\item[$E_6$] the split form $EI$ --- exactly one type  with $\frak
g_0=\frak s\frak o (5,5)\oplus \Bbb R$ and $\frak g_{-1}=\Bbb R^{16}$.
\item[$E_6$] the real form $EIV$ --- exactly one type  with
$\frak g_0=\frak s\frak o(1,9)\oplus \Bbb R$ and $\frak g_{-1}=\Bbb R^{16}$.
\item[$E_7$] the split form $EV$ --- exactly one type with $\frak g_0=
EI\oplus \Bbb R$ and $\frak g_{-1}=\Bbb R^{27}$.
\item[$E_7$] the real form $EVII$ -- exactly one type with $\frak
g_0=EIV\oplus\Bbb R$ and $\frak g_{-1}=\Bbb R^{27}$.
\end{itemize}

\subsection{Example}\label{2.10}
The {\em parabolic contact geometries} form
another important class. They correspond to $|2|$--graded Lie algebras
$\frak g$ with one--dimensional top components $\frak g_2$. Thus the regular
infinitesimal structures are equivalent to contact geometric structures,
together with the reduction of the graded tangent space to the subgroup
$G_0$ in the group of contact transformations. The only exceptions are the
so called {\em projective contact structures\/} 
($C_\ell$ series of algebras) where
more structure has to be added, see e.g. \cite{CS}. The general classification
scheme allows a simple formulation for the contact cases: The dimension one
condition on $\frak g_2$ yields the prescription which simple roots have to
be crossed while the prescribed length two of the grading gives further
restrictions. The outcome may be expressed as (see \cite{Kan,Yam}):

\begin{prop*} Each non--compact real simple Lie algebra $\frak g$ 
admits a unique grading of
contact type (up to conjugacy classes), except $\frak g$ is one of $\frak
s\frak l(2,\Bbb R)$, $\frak s\frak l(\ell,\Bbb H)$, 
$\frak s\frak p(p,q)$, $\frak s\frak o(1,q)$, 
$EIV$, $FII$ and in these cases no such gradings exist.
\end{prop*}

The  best known examples are the non--degenerate hypersurface type
{\em CR geometries} 
(with signature $(p,q)$ of the Levi form) which are exactly
the torsion free regular parabolic geometries with $\frak g=\frak s\frak
u(p+1,q+1)$, see e.g. \cite{CS}, Section 4.14--4.16. 
The real split forms of the same complex
algebras give rise to the so called {\em almost Lagrangian contact
geometries}, cf. \cite{Tak}.

\subsection{Example}\label{2.11}
The previous two lists of geometries include those with most simple
infinitesimal flag structures. The other extreme is provided by the real
parabolic geometries with most complicated flags in each tangent space,
i.e. those corresponding to the Borel subgroups $P\subset G$. Here we need
to cross out all nodes in the Satake diagram and so there must not be any black
ones. Thus all real
split forms, $\frak s\frak u(p,p)$, $\frak s\frak o(\ell-1,\ell+1)$, and EII
list all real forms which admit the right grading.

\subsection{Natural bundles}\label{2.12}
Consider a fixed parabolic geometry $(\cg, \om)$ over a manifold $M$. Then
each $P$--module $\Bbb V$ defines the associated bundle $VM=\cg\x_P \Bbb V$
over $M$. In fact, this is a functorial construction which may be 
restricted to all subcategories of parabolic geometries mentioned in \ref{2.6}.

Similarly, we may treat bundles associated to any action $P\to
\operatorname{Diff}(\Bbb S)$ on a manifold $\Bbb S$, the standard fiber for
$SM=\cg\x_P \Bbb S$. We shall meet only natural vector bundles defined by
$P$--modules in this paper, however. 

There is a special class of natural (vector) bundles defined by $G$--modules
$\Bbb W$. Such natural bundles are called {\em tractor bundles}, see
\cite{BEG, CG} for historical remarks. We shall distinguish them by the
script letters here and often omit the base manifold $M$ from the notation.
We may view each such tractor bundle $\Cal W M$ 
as associated to the extended principal fiber bundle $\tilde \cg = \cg\x_P
G$, i.e. $\Cal W =\tilde \cg\x_G\Bbb W$. Now, the Cartan connection $\om$ on
$\cg$ extends uniquely to a principal connection form $\tilde \om$ on
$\tilde \cg$, and so there is the induced linear connection on 
each such $\Cal W$.
With some more careful arguments, this construction may be extended to all
$(\frak g,P)$--modules $\Bbb W$, i.e. $P$--modules with a fixed extension of
the induced representation of $\frak p$ to a representation of $\frak g$
compatible with the $P$--action, see \cite{CG}, Section 2. One of the 
achievements of the latter paper is the equivalent treatment of the regular
parabolic geometries entirely within the framework of the tractor bundles, 
inclusive the discussion of the canonical connections.

\subsection{Adjoint tractors}\label{2.13}
It seems that the most important natural bundle is the {\em adjoint tractor
bundle} $\Cal A = \cg \x_P\frak g$ with respect to the adjoint action
$\operatorname{Ad}$ of $G$ on $\frak g$. The $P$--submodules $\frak
g^j\subset \frak g$ give rise to the filtration 
$$
\Cal A = \Cal A^{-k}\supset \Cal A^{-k+1}\supset\dots\supset \Cal A^0\supset \Cal
A^1\supset\dots\supset \Cal A^k
$$
by the natural subbundles $\ca^j=\cg\x_P\frak g^j$. Moreover, the associated
graded natural bundle (often denoted by 
the abuse of notation by the same symbol again)
$$
\operatorname{Gr}\ca=
\ca_{-k}\oplus\dots\oplus\ca_{-1}\oplus\ca_0\oplus\ca_1\oplus\dots\oplus\ca_k
$$ 
with $\ca_{j}=\ca^j/\ca^{j+1}$ is available. By the very definition, there
is the algebraic bracket on $\ca$ 
defined by means of the Lie bracket in $\frak g$ (since the Lie bracket is
$\operatorname{Ad}$-equivariant), which shows up on the graded 
bundle as
$$
\{~,~\}: \ca_i\x\ca_j\to \ca_{i+j}
.$$
For the same reason, the Killing form
defines a pairing on $\operatorname{Gr}\ca$ such that $\ca_i^*=\ca_{-i}$,
and the algebraic
codifferential $\partial^*$, cf. \ref{2.5}, defines natural algebraic
mappings 
$$
\partial^*:
\wedge^{k+1}\ca^1\otimes\operatorname\ca
\to \wedge^{k}\ca^1\otimes\operatorname\ca
$$
which are homogeneous of degree zero with respect to the gradings in
$\operatorname{Gr}\ca$.

Similarly to the notation for $\frak g$, we also write $\ca_+=\ca^1$,
$\ca_-=\ca/\ca^0$ for bundles associated either to $\cg$ or $\cg_0$.
Thus $\ca=\ca_-+\ca_0+\ca_+$, understood either as composition series induced
by the filtration, or direct sum of invariant subbundles, respectively.

\subsection{Tangent and cotangent bundles}\label{2.14}
For each parabolic geometry $(\cg,\om)$, $p:\cg\to M$, 
the absolute parallelism defines the identification
$$
\cg\x_P (\frak g/\frak p)\simeq TM,\quad \cg\x \frak g_-\ni (u,X)\mapsto
Tp(\om^{-1}(X)(u)) 
.$$
In other words, the tangent spaces $TM$ are natural bundles 
equipped with the filtrations which correspond to the Lie algebras
$\frak g_-$ viewed as the $P$--modules $\frak g/\frak p$ with  the induced
$\operatorname{Ad}$--actions. Equivalently, the tangent spaces are the
quotients
$$
TM= \ca /\ca^0
$$
of the adjoint tractor bundles. Therefore, the induced graded tangent spaces 
$\operatorname{Gr}TM$ are exactly the negative parts of the graded adjoint
tractor bundles
$$
\operatorname{Gr}TM=\ca_{-k}\oplus \dots \oplus \ca_{-1}
.$$ 
Moreover, the definition of the algebraic bracket on $\ca$ implies
immediately that the bracket induced by the Lie bracket of vector
fields on $\operatorname{Gr}TM$ for regular infinitesimal flag structures on
$M$ coincides with $\{\ ,\ \}$.

Now, the cotangent bundles clearly correspond to 
$$
T^*M=\cg\x_P \frak p_+ \simeq \ca^1
$$ 
and so the graded cotangent space is identified with 
$$
\operatorname{Gr}T^*M=\ca_1\oplus\dots\oplus\ca_k
.$$
Finally, the pairing of a one--form and a vector field is given exactly by the
canonical pairing of $\ca/\ca^1$ and $\ca^1$ induced by the 
Killing form.

\subsection{}\label{2.15}
The first important observation about the adjoint tractors and their links
to tangent and cotangent spaces is that the
curvature $K$ of the parabolic geometry $(\cg,\om)$ is in fact a section of
$\La^2(\ca/\ca^0)^*\otimes\ca$ whose 
frame form is the curvature function $\ka$. Thus, the curvature is a
two--form on the underlying manifold $M$ valued in the adjoint tractors and
all the conditions on the curvature discussed in \ref{2.6} are expressed by
natural algebraic operations on the adjoint tractors. 

The remarkable relation of both tangent and cotangent spaces to the
positive and negative parts of the adjoint tractors is the most important
tool in what follows. In particular, let us notice already here that once we
are given a reduction of the structure group $P$ of $\cg$ to its reductive
part $G_0$, the adjoint tractor bundles are identified with their graded
versions and both tangent and cotangent bundles are embedded inside of $\ca$.

\section{Weyl--structures}\label{3}

\subsection{Definition}\label{3.1}
Let $(p:\Cal G\to M,\om)$ be a parabolic geometry on a smooth manifold
$M$, and consider the underlying principal $G_0$--bundle $p_0:\Cal
G_0\to M$ and the canonical projection $\pi:\Cal G\to \Cal G_0$. A
{\em Weyl--structure\/} for $(\Cal G,\om)$ is a global
$G_0$--equivariant smooth section $\si:\Cal G_0\to\Cal G$ of $\pi$.

\begin{prop}\label{3.2}
For any parabolic geometry $(p:\Cal G\to M,\om)$, there exists a
Weyl--structure. Moreover, if $\si$ and $\hat\si$ are two
Weyl--structures, then there is a unique smooth section
$\Up=(\Up_1,\dots,\Up_k)$ of $\Cal A_1\oplus\dots\oplus\Cal A_k$ such
that 
$$
\hat\si(u)=\si(u)\exp(\Up_1(u))\dots\exp(\Up_k(u)).$$ 
Finally, each Weyl-structure $\si$ and section $\Up$ define 
another Weyl-structure $\hat\si$ by the above formula.
\end{prop}
\begin{proof}
We can choose a finite open covering $\{U_1,\dots,U_N\}$ of $M$ such
that both $\Cal G$ and $\Cal G_0$ are trivial over each $U_i$. Since
by Proposition \ref{2.3} $P$ is the semidirect product of $G_0$ and
$P_+$ it follows immediately that there are smooth $G_0$--equivariant 
sections $\si_i:p_0^{-1}(U_i)\to p^{-1}(U_i)$. Moreover, we can find
open subsets $V_i$ such that $\bar V_i\subset U_i$ and such that
$\{V_1,\dots,V_N\}$ still is a covering of $M$. 

Now from Proposition \ref{2.3} and the Baker--Campbell--Hausdorff
formula it follows that there is a smooth mapping
$\Ps:p_0^{-1}(U_1\cap U_2)\to\frak p_+$ such that
$\si_2(u)=\si_1(u)\exp(\Ps(u))$. Equivariance of $\si_1$ and $\si_2$
immediately implies that $\Ps(u\cdot g)=\Ad(g^{-1})(\Ps(u))$ for all
$g\in G_0$. Now let $f:M\to [0,1]$ be a smooth function with support
contained in $U_2$, which is identically one on $V_2$ and define
$\si:p_0^{-1}(U_1\cup V_2)\to p^{-1}(U_1\cup V_2)$ by
$\si(u)=\si_1(u)\exp(f(p_0(u))\Ps(u))$ for $u\in U_1$ and by
$\si(u)=\si_2(u)$ for $u\in V_2$. Then obviously these two definitions
coincide on $U_1\cap V_2$, so $\si$ is smooth. Moreover, from the
equivariance of the $\si_i$ and of $\Ps$ one immediately concludes
that $\si$ is equivariant. Similarly, one extends the section next to
$U_1\cup V_2\cup V_3$ and by induction one reaches a globally defined
smooth equivariant section. 

If $\hat\si$ and $\si$ are two global equivariant sections, then
applying Proposition \ref{2.3} directly, we see that there are smooth
maps $\Up_i:\Cal G_0\to\frak g_i$ for $i=1,\dots,k$ such that
$\hat\si(u)=\si(u)\exp(\Up_1(u))\dots\exp(\Up_k(u))$. As above,
equivariance of $\hat\si$ and $\si$ implies that $\Up_i(u\cdot
g)=\Ad(g^{-1})(\Up_i(u))$ for all $g\in G_0$. 
Hence, $\Up_i$ corresponds to a smooth
section of $\Cal A_i$. The last statement of the Proposition is obvious now. 
\end{proof}

\subsection{Weyl connections}\label{3.3}
We can easily relate a Weyl--structure $\si:\Cal G_0\to\Cal G$ to
objects defined on the manifold $M$ by considering the pullback
$\si^*\om$ of the Cartan connection $\om$ along the section
$\si$. Clearly, $\si^*\om$ is a $\frak g$--valued one--form on $\Cal
G_0$, which by construction is $G_0$--equivariant,
i.e. $(r^g)^*(\si^*\om)=\Ad(g^{-1})\o\si^*\om$ for all $g\in
G_0$. Since $\Ad(g^{-1})$ preserves the grading of $\frak g$, in fact
each component $\si^*\om_i$ of $\si^*\om$ is a $G_0$--equivariant one
form with values in $\frak g_i$. 

Now consider a vertical tangent vector on $\Cal G_0$, i.e.~the value
$\ze_A(u)$ of a fundamental vector field corresponding to some
$A\in\frak g_0$. Since $\si$ is $G_0$--equivariant, we conclude that
$T_u\si\cdot\ze_A(u)=\ze_A(\si(u))$, where the second fundamental
vector field is on $\Cal G$. Consequently, we have
$\si^*\om(\ze_A)=\om(\ze_A)=A\in\frak g_0$. Thus, for
$i\neq 0$ the form $\si^*\om_i$ is horizontal, while $\si^*\om_0$
reproduces the generators of fundamental vector fields. 

From this observation, it follows immediately, that for $i\neq 0$, the
form $\si^*\om_i$ descends to a smooth one form on $M$ with values in
$\Cal A_i$, which we denote by the same symbol, while $\si^*\om_0$
defines a principal connection on the bundle $\Cal G_0$. This connection is
called the {\em Weyl connection} of the Weyl structure $\si$. 

\subsection{Soldering forms and Rho-tensors}\label{3.4}
We view the
positive components of $\si^*\om$ as a one--form 
$$
\Rho=\si^*(\om_+)\in\Om^1(M;\ca_1\oplus\dots\oplus\ca_k)
$$ 
with values in the bundle $\Cal A_1\oplus\dots\oplus\Cal A_k$.
We call it the {\em
Rho--tensor\/} of the Weyl--structure $\si$. This is a generalization
of the tensor $\Rho_{ab}$ well known in conformal geometry. 

Since $\om$ restricts to a linear isomorphism in each tangent space of
$\Cal G$, we see that the form 
$$
\si^*\om_-=(\si^*\om_{-k},\dots,\si^*\om_{-1})\in\Om^1(M,\Cal
A_{-k}\oplus\dots\oplus\Cal A_{-1})
$$
induces an isomorphism 
$$
TM\cong \Cal A_{-k}\oplus\dots\oplus\Cal
A_{-1}\cong \operatorname{Gr}TM.
$$ 
We will denote this isomorphism by 
$$
\xi\mapsto
(\xi_{-k},\dots,\xi_{-1})\in \ca_{-k}\oplus\dots\oplus\ca_{-1}
$$ 
for $\xi\in TM$.  In particular, each fixed $u\in\cg_0$ provides the
identification of $T_{p_0(u)}M\cong\frak g_-$ compatible with the grading.
Thus, the choice of a Weyl structure $\si$ provides a reduction of the 
structure group of $TM$ to $G_0$ (by means of the soldering form 
$\si^*\om_-$ on $\cg_0$), the linear connection on $M$ 
(the Weyl connection $\si^*\om_0$),
and the Rho--tensor $\Rho$.

\subsection{Remarks}\label{3.5}
As discussed in \ref{2.7}--\ref{2.8} above, 
there is the underlying frame form of length
one on $\Cal G_0$ which is the basic structure from which the whole parabolic
geometry $(\cg,\om)$ may be reconstructed, with exceptions mentioned
explicitly in \ref{2.9} and \ref{2.10}. 
By definition, for $i<0$ and $\xi\in T^i\Cal G_0$
this frame form can be computed by choosing any lift of $\xi$ to a
tangent vector on $\Cal G$ and then taking the $\frak g_i$--component
of the value of $\om$ on this lift. In particular, we can use
$T\si\cdot\xi$ as the lift, which implies that the restriction of
$\si^*\om_i$ (viewed as a form on $\Cal G_0$) to $T^i\Cal G_0$
coincides with the $\frak g_i$--component of the frame form of length
one. This in turn implies that the restriction of $\si^*\om_i$ (viewed
as a form on $M$) to $T^iM$ coincides with the canonical projection
$T^iM\to\Cal A_i=T^iM/T^{i+1}M$.

There is also another interpretation of the objects on $M$ induced by the
choice of a Weyl--structure that will be very useful in the
sequel. Namely, consider the form 
$$
\si^*\om_{\leq 0}=\si^*\om_{-k}\oplus\dots\oplus\si^*\om_0
\in\Om^1(\Cal G_0,\frak g_{-k}\oplus\dots\oplus\frak g_0).
$$ 
We have seen above that this form is $G_0$--equivariant, it reproduces
the generators of fundamental vector fields, and restricted to each
tangent space, it is a linear isomorphism. Thus $\si^*\om_{\leq 0}$
defines a Cartan connection on the principal $G_0$--bundle $p_0:\Cal
G_0\to M$. In the case of the irreducible parabolic geometries, these
connections are classical affine connections on the tangent space $TM$ 
belonging to its reduced structure group $G_0$.

\subsection{Bundles of scales}\label{3.6}
As we have seen in \ref{3.3}, \ref{3.4} above, choosing a Weyl--structure
$\si:\Cal G_0\to\Cal G$ leads to several objects on the manifold
$M$. Now the next step is to show that in fact a small part of these
data is sufficient to completely fix the Weyl--structure. More
precisely, we shall see below that even the linear connections induced
by the Weyl connection $\si^*\om_0$ on certain oriented
line bundles suffice to pin down the Weyl--structure. Equivalently,
one can use the corresponding frame bundles, which are principal
bundles with structure group $\Bbb R^+$. The principal bundles
appropriate for this purpose are called {\em bundles of scales\/}.

To define these bundles, we have to make a few observations: A
principal $\Bbb R^+$--bundle associated to $\Cal G_0$ is determined by
a homomorphism $\la:G_0\to \Bbb R^+$. The derivative of this
homomorphism is a linear map $\la':\frak g_0\to\Bbb R$. Now $\frak
g_0$ splits as the direct sum $\frak z(\frak g_0)\oplus \frak
g_0^{ss}$ of its center and its semisimple part, and $\la'$
automatically vanishes on the semisimple part. Moreover, as discussed in
\ref{2.2} the restriction of the Killing form $B$ of
$\frak g$ to the subalgebra $\frak g_0$ is non--degenerate, and one
easily verifies that this restriction respects the above splitting. In
particular, the restriction of $B$ to $\frak z(\frak g_0)$ is still
non--degenerate and thus there is a unique element $E_\la\in\frak
z(\frak g_0)$ such that $\la'(A)=B(E_\la,A)$ for all $A\in\frak g_0$. 

Next, the action of the element $E_\la\in\frak z(\frak g_0)$ on any
$G_0$--irreducible representation commutes with the action of $G_0$,
and thus is given by a scalar multiple of the identity by Schur's
lemma. 

\subsection*{Definition}
An element $E_\la$ of $\frak z(\frak g_0)$ is called a {\em scaling
element\/} if and only if $E_\la$ acts by a nonzero real scalar on
each $G_0$--irreducible component of $\frak p_+$. 
A {\em bundle of scales\/} is a principal $\Bbb R^+$ bundle $\Cal
L^\la\to M$ which is associated to $\Cal G_0$ via a homomorphism
$\la:G_0\to\Bbb R^+$, whose derivative is given by
$\la'(A)=B(E_\la,A)$ for some scaling element $E_\la\in\frak z(\frak
g_0)$.

Having given a fixed choice of a bundle $\Cal L^\la$ of scales, a
{\em (local) scale\/} on $M$ is a (local) smooth section of $\Cal
L^\la$. 

\begin{prop}\label{3.7}
Let $G$ be a fixed semisimple Lie group, whose Lie algebra $\frak g$
is endowed with a $|k|$--grading. Then the following holds:\newline
(1)\ There are scaling elements in $\frak z(\frak g_0)$.\newline
(2)\ Any scaling element $E_\la\in\frak z(\frak g_0)$ gives rise to a
canonical bundle $\Cal L^\la$ of scales over each manifold endowed
with a parabolic geometry of the given type.\newline
(3)\ Any bundle of scales admits global smooth sections, i.e.~there
always exist global scales.
\end{prop}
\begin{proof}
(1) The grading element $E\in\frak z(\frak g_0)$, cf. \ref{2.2}, 
    acts on $\frak
    g_i$ by multiplication with $i$, so it is a scaling element. More
    generally, one can consider the subspace of $\frak z(\frak g_0)$
    of all elements which act by real scalars on each irreducible
    component of $\frak p_+$. Then each irreducible component
    determines a real valued functional and thus a hyperplane in
    that space, and the complement of these finitely many hyperplanes
    (which is open and dense) consists entirely of scaling
    elements.\newline
(2) Let $\frak p_+=\oplus\frak p^\al$ be the decomposition of $\frak
    p_+$ into $G_0$--irreducible components, and for a fixed grading
    element $E_\la$ denote by $a_\al$ the scalar by which $E_\la$ acts
    on $\frak p^\al$. The adjoint action defines a smooth homomorphism
    $G_0\to \prod_{\al}GL(\frak p^\al)$, whose components 
    we write as $g\mapsto
    \Ad^\al(g)$. Then consider the homomorphism $\la:G_0\to\Bbb R^+$
    defined by
    $$
\la(g):=\prod_{\al}|\det(\Ad^\al(g))|^{2a_\al}.
    $$ 
    The derivative of this homomorphism is given by
    $\la'(A)=\sum_\al 2a_\al
    \text{tr}(\ad(A)|_{\frak p^\al})$. Now $\frak
    g_{-}=\oplus_\al(\frak p^\al)^*$, and $E_\la$ acts on $(\frak
    p^\al)^*$ by $-a_\al$ and on $\frak g_0$ by zero, and thus
    $B(E_\la,A) = \operatorname{tr}(\ad(A)\o\ad(E_\la))=\sum_\al
a_\al\operatorname{tr}(\ad(A)|_{\frak p^\al}) - \sum_\al
a_\al\operatorname{tr}(\ad(A)|_{(\frak p^\al)^*}) = \la'(A)$.
\newline
(3) This is just due to the fact that orientable real line bundles and
    thus principal $\Bbb R^+$--bundles are automatically trivial and
    hence admit global smooth sections. 
\end{proof}

\begin{lem}\label{3.8}
Let $\si:\Cal G_0\to\Cal G$ be a Weyl--structure for parabolic geometry 
$(\cg\to M,\om)$ and let $\Cal L^\la$ be a bundle of scales.\newline
(1) The Weyl connection $\si^*\om_0\in\Om_1(\Cal G_0,\frak g_0)$ 
induces a principal connection on the bundle of scales $\Cal L^\la$. \newline
(2) $\Cal L^\la$ is naturally identified with $\Cal G_0/\ker(\la)$, 
the orbit space of
the free right action of the normal subgroup $\ker(\la)\subset G_0$ on
$\Cal G_0$. \newline
(3) The form $\la'\o\si^*\om_0\in\Om^1(\Cal G_0)$ descends to 
the connection form of the induced
principal connection on $\Cal L^\la=\cg_0/\ker(\la)$. \newline
(4) The composition of $\la'$ with the curvature form of $\si^*\om_0$
descends to the curvature of the induced connection on $\Cal L^\la$.
\end{lem}
\begin{proof}
All claims are straightforward consequences of the definitions. 
\end{proof}

To see that the Weyl--structure $\si$ is actually uniquely determined by the
induced principal connection on $\Cal L^\la$ (cf. Theorem \ref{3.12} below),
we have to compute how the principal connection $\si^*\om_0$ changes when we
change $\si$. For later use, we also compute how the other objects induced
by $\si$ change under the change of the Weyl--structures. So let us assume
that $\hat\si$ is another Weyl--structure and
$\Up=(\Up_1,\dots,\Up_k)$ is the section of $\Cal
A_1\oplus\dots\oplus\Cal A_k$ characterized by
$\hat\si(u)=\si(u)\exp(\Up_1(u))\dots\exp(\Up_k(u))$.

We shall use the convention that we simply denote quantities corresponding
to $\hat\si$ by hatted symbols and quantities corresponding to $\si$ by
unhatted symbols. Consequently, $(\xi_{-k},\dots,\xi_{-1})$ and
$(\hat\xi_{-k},\dots,\hat\xi_{-1})$ denote the splitting of $\xi\in
TM$ according to $\si$, respectively $\hat\si$, and $\Rho$ and $\hat\Rho$ are
the Rho--tensors. Finally, let us consider any vector bundle $E$
associated to the principal bundle $\Cal G_0$. Then for any Weyl--structure
the corresponding principal connection on $\Cal G_0$ induces a linear
connection on $E$, which is denoted by $\nabla$ for $\si$ and by
$\hat\nabla$ for $\hat\si$.

To write the formulae efficiently, we need some further notation. By
$\underline{j}$ we denote a sequence $(j_1,\dots,j_k)$ of nonnegative
integers, and we put
$\|\underline{j}\|=j_1+2j_2+\dots+kj_k$. Moreover, we define
$\underline{j}!=j_1!\dots j_k!$ and
$(-1)^{\underline{j}}=(-1)^{j_1+\dots+j_k}$,  
and we define $(\underline{j})_m$ to
be the subsequence $(j_1,\dots ,j_m)$ of $\underline{j}$. By $0$ we
denote sequences of any length consisting entirely of zeros.

\begin{prop}\label{3.9}
Let $\si$ and $\hat\si$ be two Weyl--structures related by
$$
\hat\si(u)=\si(u)\exp(\Up_1(u))\dots\exp(\Up_k(u)),
$$
where $\Up=(\Up_1,\dots,\Up_k)$ 
is a smooth section of $\Cal A_1\oplus\dots\oplus\Cal
A_k$. Then we have: 
\begin{align}\label{splittransf}
\hat\xi_i&=\sum_{
\|\underline{j}\|+\ell=i}
\frac{(-1)^{\underline{j}}}{\underline{j}!}\ad(\Up_k)^{j_k}\o\dots
\o\ad(\Up_1)^{j_1}(\xi_\ell),\\
\hat\Rho_i(\xi)&=\sum_{
\|\underline{j}\|+\ell=i}
\frac{(-1)^{\underline{j}}}
{\underline{j}!}\ad(\Up_k)^{j_k}\o\dots\o\ad(\Up_1)^{j_1}(\xi_\ell)+
\label{Rhotransf}\\
&\quad\sum_{
\|\underline{j}\|+\ell=i}
\frac{(-1)^{\underline{j}}}{\underline{j}!}
\ad(\Up_k)^{j_k}\o\dots\o\ad(\Up_1)^{j_1}(\Rho_\ell(\xi))+\notag\\
&\quad\sum_{m=1}^k\sum_{\substack{
(\underline{j})_{m-1}=0\\ 
m+\|\underline{j}\|=i}}
\frac{(-1)^{\underline{j}}}{(j_m+1)\underline{j}!}
\ad(\Up_k)^{j_k}\o\dots\o\ad(\Up_m)^{j_m}(\nabla_\xi\Up_m),
\notag\\
\noalign{\noindent 
where $\ad$ denotes the adjoint action with respect to the algebraic
bracket $\{\ ,\ \}$. \newline\indent
If $E$ is an associated vector bundle to the principal bundle $\Cal
G_0$, then we have: }
\hat\nabla_\xi s&=\nabla_\xi
s+\sum_{
\|\underline{j}\|+\ell=0} 
\frac{(-1)^{\underline{j}}}{\underline{j}!}(\ad(\Up_k)^{j_k}\o\dots\o
\ad(\Up_1)^{j_1}(\xi_\ell))\bullet s,
\label{conntransf}
\end{align}
where $\bullet$ denotes the map $\Cal A_0\x E\to E$ induced by the
action of $\frak g_0$ on the standard fiber of $E$.
\end{prop}
\begin{proof}
The essential part of the proof is to compute the tangent map
$T_u\hat\si$ in a point $u\in\Cal G_0$. By definition,
$\hat\si(u)=\si(u)\exp(\Up_1(u))\dots\exp(\Up_k(u))$. Thus, we can
write the evaluation of the tangent map, $T_u\hat\si\cdot\xi$, as the sum of
$T_{\si(u)}r^g\cdot T_u\si\cdot\xi$, where
$g=\exp(\Up_1(u))\dots\exp(\Up_k(u))\in P_+$, and the derivative at
$t=0$ of 
$$
\si(u)\exp(\Up_1(c(t)))\dots\exp(\Up_k(c(t))),
$$ 
where $c:\Bbb R\to \cg_0$ is a smooth curve with $c(0)=u$ and
$c'(0)=\xi$. By construction, the latter derivative lies in the kernel
of $T\pi$, where $\pi:\Cal G\to\Cal G_0$ is the projection, so we can
write it as $\ze_{\Ph(\xi)}(\hat\si(u))$ for suitable $\Ph(\xi)\in\frak p_+$.

Now, for $\xi\in T_u\Cal G_0$, we have
$\hat\si^*\om(\xi)=\om(\hat\si(u))(T_u\hat\si\cdot\xi)$. By
equivariance of the Cartan connection $\om$, we get $\om(\si(u)\cdot
g)(Tr^g\cdot T\si\cdot
\xi)=\Ad(g^{-1})(\om(u)(T\si\cdot\xi))$. Consequently,
$$
\hat\si^*\om(\xi)=\Ad(g^{-1})(\si^*\om(\xi))+\Ph(\xi). 
$$
Since $\Ph(\xi)\in\frak p_+$, this term affects only the transformation of
the Rho--tensor, and does not influence the changes of $\si^*\om_i$
for $i\leq 0$. In particular, for the components $\hat\si^*\om_i$ with
$i<0$, we only have to take the part of the right degree in 
\begin{equation}
e^{\ad(-\Up_k(u))}\o\dots\o e^{\ad(-\Up_1(u))}(\si^*\om(u)(\xi)),
\label{3.9.4}\end{equation}
and expanding the exponentials, this immediately leads to formula
\eqref{splittransf}. 

To compute the change in the connection, we have to notice that
$\hat\si^*\om_0(\xi)$ is the component of degree zero in \eqref{3.9.4}
above. Consequently, if we apply $\hat\si^*\om_0$ to the horizontal
lift of a tangent vector on $M$, the outcome is just this degree zero
part. Otherwise put, the horizontal lift with respect to
$\hat\si^*\om_0$ is obtained by subtracting the fundamental
vector field corresponding to the degree zero part of \eqref{3.9.4} from the
horizontal lift with respect to $\si^*\om_0$. Applying such horizontal
vector field to a
smooth $G_0$--equivariant function with values in any
$G_0$--representation and taking into account that a fundamental
vector fields acts on such functions by the negative of its generator acting 
on the values, this immediately leads to formula \eqref{conntransf}
by expanding the exponentials.

Finally, we have to deal with the change of the Rho--tensor. Recall
that we view this as a tensor on the manifold $M$, so we can compute
$\hat\Rho_i(\xi)$ by applying $\hat\si^*\om_i$ to any lift of
$\xi$. In particular, we may use the horizontal lift $\xi^h$ with
respect to $\si^*\om_0$, so we may assume $\si^*\om_0(\xi)=0$. But
then expanding the exponentials in \eqref{3.9.4} and taking the part of degree
$i$ we see that we exactly get the first two summands in formula
\eqref{Rhotransf}. Thus we are left with proving that the last summand
corresponds to $\Ph(\xi)$. For this aim, let us rewrite the curve that we
have to differentiate as
$$
\hat\si(u)\exp(-\Up_k(u))\dots\exp(-\Up_1(u))
\exp(\Up_1(c(t)))\dots\exp(\Up_k(c(t))). 
$$
Differentiating this using the product rule we get a sum of terms in
which one $\Up_i$ is differentiated, while all others have to be
evaluated at $t=0$, i.e.~in $u$. So each of these terms reads as the
derivative at $t=0$ of 
$$
\hat\si(u)\cdot\text{conj}_{\exp(-\Up_k(u))}\o\dots\o
\text{conj}_{\exp(-\Up_{i+1}(u))}\bigl(\exp(-\Up_i(u))\exp(\Up_i(c(t)))\bigr),
$$
where $\text{conj}_g$ denotes the conjugation by $g$, i.e.~the map
$h\mapsto ghg^{-1}$. This expression is just the principal right
action by the value of a smooth curve in $P$ which maps zero to the
unit element, so its result is exactly the value at $\hat\si(u)$ of 
the fundamental vector field generated by the derivative at zero of
this curve. This derivative is clearly obtained by applying
$$
e^{\ad(-\Up_k(u))}\o\dots\o e^{\ad(-\Up_{i+1}(u))}
$$
to the derivative at zero of $t\mapsto
\exp(-\Up_i(u))\exp(\Up_i(c(t)))$. By \cite{KMS}, 4.26, and the chain
rule, the latter derivative equals the left logarithmic
derivative of $\exp$ applied to the derivative at zero of
$t\mapsto \Up_i(c(t))$. Moreover, the proof of \cite{KMS}, Lemma 4.27,
can be easily adapted to the left logarithmic derivative, showing that
this gives  
$$
\sum_{p=0}^\infty\frac{(-1)^p}{(p+1)!}\ad(\Up_i(u))^p(\xi^h\cdot\Up_i).  
$$
Finally, we have to observe that $\xi^h\cdot\Up_i$ corresponds to
$\nabla_\xi\Up_i$ and to sort out the terms of the right degree in order 
to get the remaining summand in \eqref{Rhotransf}. 
\end{proof}

\subsection{Example}\label{3.10} 
For all irreducible parabolic geometries, the formulae
from Proposition \ref{3.9} become extremely simple. In fact they coincide
completely with the known ones in the conformal Riemannian geometry: The
grading of $TM$ is trivial, the connection transforms as
$$
\hat\nabla_\xi s= \nabla_\xi s - \{\Up,\xi\}\bullet s,
$$
where $\Up$ is a section of $\ca_1=T^*M$, and the bracket of $\Up$ and $\xi$ is
a field of endomorphisms of $TM$ acting on $s$ in an obvious way. Indeed,
there are no more terms on the right--hand side of
\ref{3.9}\eqref{conntransf} which make sense.
Next, the Rho--tensor transforms as 
$$
\hat\Rho(\xi) = \Rho(\xi) + \nabla_\xi\Up + \tfrac12\{\Up,\{\Up,\xi\}\} 
.$$
The formulae for the $|2|$--graded examples are a bit more complicated. The
splitting of $TM$ and the connection and Rho--tensors change as follows
\begin{align*}
\hat\xi_{-2} =\ &\xi_{-2}
\\
\hat\xi_{-1} =\ &\xi_{-1} - \{\Up_1,\xi_{-2}\}
\\
\hat\nabla_\xi s =\ & \nabla_\xi s + \bigl(
\tfrac12\{\Up_1,\{\Up_1,\xi_{-2}\}\}-\{\Up_2,\xi_{-2}\}  - 
\{\Up_1,\xi_{-1}\})\bullet s,
\\
\hat\Rho_1(\xi) =\ &\Rho_1(\xi) -\tfrac16 \{\Up_1,\{\Up_1,\{\Up_1,\xi_{-2}\}\}\}
+ \{\Up_2,\{\Up_1,\xi_{-2}\}\} +
\\
&\tfrac12\{\Up_1,\{\Up_1,\xi_{-1}\}\} -
\{\Up_2,\xi_{-1}\} +\nabla_{\xi}\Up_1 
\\
\hat\Rho_2(\xi) =\ &\Rho_2(\xi) -\{\Up_1,\Rho_1(\xi)\} + \nabla_\xi \Up_2
-\tfrac12\{\Up_1,\nabla_\xi\Up_1\} +
\\
&\tfrac1{24}\ad(\Up_1)^4(\xi_{-2}) -\tfrac12 \{\Up_2,\{\Up_1,
\{\Up_1,\xi_{-2}\}\}\}
+ \tfrac12\{\Up_2,\{\Up_2,\xi_{-2}\}\} -
\\
& \tfrac16\ad(\Up_1)^3(\xi_{-1}) + \{\Up_2,\{\Up_1,\xi_{-1}\}\}
.\end{align*}

\subsection{Remark}\label{3.11} 
In applications, one is often interested in questions about the dependence
of some objects on the choice of the Weyl--structures and then the 
infinitesimal form of the available change of the splittings, Rho's and
connections is important. In our terms, this amounts to sorting
out the terms in formulae \ref{3.9}\eqref{splittransf}--\eqref{conntransf}
which are linear in upsilons. Thus,
the infinitesimal version of Proposition \ref{3.9} for the variations
$\delta\xi_i$, $\delta\nabla$, and $\delta\Rho_i$ reads 
\begin{align}\label{splittransfinf}
\delta\xi_i =\ &- \{\Up_1,\xi_{i-1}\} - \dots - \{\Up_{k+i},\xi_{-k}\}
\\
\delta\Rho_i(\xi)=\ & \nabla_\xi\Up_i
- \{\Up_1,\Rho_{i-1}(\xi)\} - \dots -
\{\Up_{i-1},\Rho_1(\xi)\}-
\label{Rhotransfinf}\\
&\{\Up_{i+1},\xi_{-1}\} - \dots -
\{\Up_{k},\xi_{-k+i}\}
\notag\\
\delta\nabla_\xi s=\ &- (\{\Up_1,\xi_{-1}\} + \dots + \{\Up_{k},\xi_{-k}\})\bullet s
\label{conntransfinf}
.\end{align}

\medskip
\subsection{}\label{3.12}
Proposition \ref{3.9} not only allows us to show that a
Weyl--structure is uniquely determined by the induced connection on
any bundle of scales, but it also leads to a description of the Cartan
bundle $p:\Cal G\to M$. To get this description, recall that for any
principal bundle $E\to M$ there is a bundle $QE\to M$ whose sections
are exactly the principal connections on $E$, see \cite{KMS}, 17.4.

\begin{thm*}
Let $p:\Cal G\to M$ be a parabolic geometry on $M$, and let $\Cal
L^\la\to M$ be a bundle of scales.\newline
(1) Each Weyl--structure $\si:\Cal G_0\to\Cal G$ determines the principal
connection on $\Cal L^\la$ induced by the Weyl connection 
$\si^*\om_0$. This defines a bijective correspondence 
between the set of Weyl--structures and the set of
principal connections on $\Cal L^\la$.\newline
(2) There is a canonical isomorphism $\Cal G\cong p_0^*Q\Cal L^\la$,
where $p_0:\Cal G_0\to M$ is the projection. 
Under this isomorphism, the choice of a Weyl structure $\si:\cg_0\to \cg$ 
is the pullback of the
principal connection on the bundle of scales $\Cal L^\la$, viewed as a
section $M\to Q\Cal L^\la$. 
Moreover, the principal action of $G_0$ is the canonical
action on $p_0^*Q\Cal L^\la$ induced from the action on $\cg_0$,
while the action of $P_+$ is described by equation
\eqref{conntransf} from Proposition \ref{3.9}.
\end{thm*}
\begin{figure}
$$
\xymatrix@R6mm{
{\text{Cartan connection $\om$}}
\ar@{.}[d]
\\
{\cg}
\ar[r]\ar[dd]
& {Q\Cal L^\la}
\ar[dd]
\\
\\
{\cg_0=\cg/P_+}
\ar@(lu,ld)[uu]^{\text{Weyl-structure $\si$}}
\ar[r]
& M
\ar@(ur,dr)[uu]_{\text{\begin{tabular}{l} principal
connection\\induced by $\si^*\om_0$\end{tabular}}}
\\
{\begin{tabular}{c} soldering form $\si^*\om_-\in\Om^1(\cg_0,\frak g_-)$
\\ Weyl connection $\si^*\om_0\in\Om^1(\cg_0,\frak g_0)$
\\ Rho--tensor $\Rho=\si^*\om_+\in\Om^1(\cg_0.\frak
p_+)$\end{tabular}}
\ar@{.}[u]
& {\begin{tabular}{c}
$\si^*\om_-\in \Om^1(M;\ca_-)$
\\ $\Rho\in \Om^1(M;\ca_+)$
\end{tabular}}
\ar@{.}[u]
}
$$
\caption{Pullback diagram with further objects related to Weyl--structures}
\end{figure}
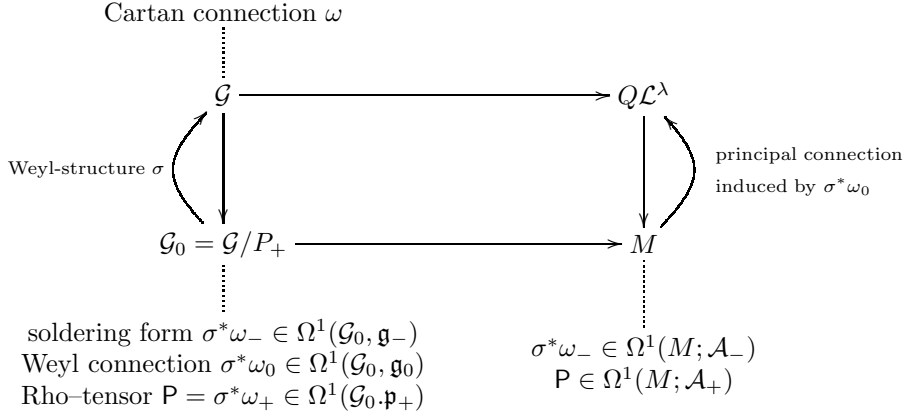

\begin{proof}
(1) Consider the map $\la':\frak g_0\to\Bbb R$ defining the bundle
   $\Cal L^\la$ of scales. Take elements $Z\in\frak p_+$ and
   $X\in\frak g_-$, and consider $\la'([Z,X])$. By assumption, this is
   given by $B(E_\la,[Z,X])=B([E_\la,Z],X)$ for some scaling element
   $E_\la\in\frak z(\frak g_0)$. Hence if we assume that $Z$ lies in a
   $G_0$--irreducible component of $\frak p_+$ this is just a nonzero
   real multiple of $B(Z,X)$. In particular, this implies that for
   each $0\neq Z\in\frak p_i$, we can find an element $X\in\frak
   g_{-i}$, such that $\la'([Z,X])\neq 0$. Moreover, since
   $E_\la\in\frak z(\frak g_0)$ we get $\Ad(g)(E_\la)=E_\la$ for all
   $g\in G_0$ and this immediately implies that mapping
   $Z\in\frak g_i$ to $X\mapsto\la'([Z,X])$ induces an isomorphism
   $\frak g_i\cong\frak g_{-i}^*$ of $G_0$--modules.

   To prove (1), we may as well use the induced linear connection on
   the line bundle $L^\la=\Cal L^\la\x_{\Bbb R^+}\Bbb R$ corresponding
   to the standard representation. For this bundle, the map $\bullet$
   from Proposition \ref{3.9} is clearly given by $(A\bullet
   s)(x)=\la'(A(x))s(x)$, where we denote by $\la':\Cal A_0\to M\x
   \Bbb R$ also the mapping induced by $\la':\frak g_0\to\Bbb R$.

   We first claim that the map from Weyl--structures to linear
   connections is injective. So assume that $\si$ and $\hat\si$ induce
   the same linear connection on $L^\la$ and let $\Up$ be the section
   of $\Cal A_1\oplus\dots\oplus\Cal A_k$ describing the change from
   $\si$ to $\hat\si$. For $\xi\in T^{-1}M$, we have $\xi_i=0$ for all
   $i<-1$, hence formula \eqref{conntransf} of \ref{3.9} reduces to
   $\hat\nabla_\xi s=\nabla_\xi s+\la'(\{\Up_1,\xi\})s$ in this
   case. If $\Up_1$ would be nonzero, then by the above argument we
   could find $\xi$ such that $\la'(\{\Up_1,\xi\})\neq 0$, which would
   contradict $\hat\nabla=\nabla$, so $\Up_1$ must be identically
   zero. But then for $\xi$ in $T^{-2}M$, the change reduces to
   $\hat\nabla_\xi s=\nabla_\xi s+\la'(\{\Up_2,\xi_{-2}\})s$ and as
   above, we conclude that $\Up_2$ is identically zero. Inductively,
   we get $\Up=0$ and thus $\hat\si=\si$.

   To see surjectivity, assume that $\hat\nabla$ is any linear
   connection on $L^\la$, and let $\si$ be any Weyl--structure with
   induced linear connection $\nabla$ on $L^\la$. Then there is a
   one--form $\tau\in\Om^1(M)$ such that $\hat\nabla_\xi s=\nabla_\xi
   s+\tau(\xi)s$. Restricting $\tau$ to $T^{-1}M$, we can find a
   unique smooth section $\Up_1$ of $\Cal A_1$ such that
   $\tau(\xi)=-\la'(\{\Up_1,\xi\})$ for all $\xi$ in $T^{-1}M$. Next,
   consider the map $T^{-2}M\to M\x\Bbb R$ given by 
$$
\xi\mapsto\tau(\xi)+\la'(\{\Up_1,\xi_{-1}\})
-\frac{1}{2}\la'(\{\Up_1,\{\Up_1,\xi_{-2}\}\}),
$$ 
   where the $\xi_i$ are the components of $\xi$ with respect to the
   Weyl--structure $\si$. By construction, this vanishes on $T^{-1}M$,
   so it factors to a map defined on $\Cal A_{-2}$, and thus there is
   a unique section $\Up_2$ of $\Cal A_2$ such that it equals
   $-\la'(\{\Up_2,\xi_{-2}\})$. Inductively, we find a section $\Up$
   such that the Weyl--structure $\hat\si$ corresponding to $\si$ and
   $\Up$ induces the linear connection $\hat\nabla$, cf. formula
\ref{3.9}\eqref{conntransf}.\newline 
(2) Consider any point $u\in\Cal G$. Proposition
   \ref{3.2} implies that there is a
   Weyl--structure $\si:\Cal G_0\to\Cal G$ such that
   $u=\si(\pi(u))$. If $\nabla$ is the linear connection on $L^\la$
   induced by $\si$, then we see from Proposition \ref{3.9} that the
   value of $\nabla_\xi s(p(u))$ for a vector field $\xi$ on $M$ and a
   section $s$ of $L^\la$ depends only on $\si(p(u))$, since its
   change under a change of the Weyl--structure depends only on the
   value of $\Up$ in $p(u)$. Thus, mapping $u$ to the value at $p(u)$
   of the principal connection on $\Cal L^\la$ induced by $\si^*\om_0$
   is independent of the choice of $\si$, so we get a well defined
   bundle map from the bundle $\Cal G\to\Cal G_0$ to the bundle $Q\Cal
   L^\la\to M$ covering the projection $p_0:\Cal G_0\to M$. Moreover,
   from part (1) of this proof it follows that this map induces
   isomorphisms in each fiber, so it leads to an isomorphism $\Cal
   G\to p_0^*Q\Cal L^\la$ of bundles over $\Cal G_0$. Obviously, the
$G_0$--equivariant sections of $p_0^*Q\Cal L^\la\to \cg_0$ correspond
exactly to the induced principal connections on $\Cal L^\la$, i.e.
the sections of $Q\Cal L^\la\to M$.
 
   In order to describe the principal action of $P$ on $p_0^*Q\Cal L^\la$
   obtained by the above isomorphism, one just has to note that for
   $u\in\Cal G_0$ and $g\in G_0$ the fibers of $p_0^*Q\Cal L^\la$ over
   $u$ and $u\cdot g$ are canonically isomorphic since
   $p_0(u)=p_0(u\cdot g)$. Thus, the principal right action of $G_0$
   is simply given by acting on $\Cal G_0$. On the other hand, fix
   $u\in\Cal G_0$ and an element $\exp(Z_1)\dots\exp(Z_k)\in P_+$ for
   $Z_i\in\frak g_i$. Via $u$, the element $Z_i$ corresponds to an
   element $\Up_i\in\Cal A_i$ at the point $p_0(u)$. Then the
   principal right action of $P_+$ is  described by the formula
   \eqref{conntransf} of Proposition \ref{3.9} as required.
\end{proof}

\subsection{Closed and exact Weyl--structures}\label{3.13}
Let us fix a bundle of scales 
$\Cal L^\la$ for some parabolic geometry. Then the
bijective correspondence between Weyl--structures and principal
connections on $\Cal L^\la$ immediately leads to two distinguished
subclasses of Weyl--structures. Namely, we call a Weyl--structure
$\si:\Cal G_0\to\Cal G$ {\em closed\/}, if the induced principal
connection on $\Cal L^\la$ (or equivalently the induced linear
connection $\nabla$ on $L^\la$) is flat.

Moreover, by Proposition \ref{3.7} the bundle $\Cal L^\la$ of scales
admits global smooth sections, and any such section gives rise to a
flat principal connection on $\Cal L^\la$ (which in addition has
trivial holonomy) and hence to a closed Weyl--structure. The closed
Weyl--structures induced by such global sections are called {\em exact}. 

Note that in the case of conformal structures, the canonical choice
for the bundle of scales is simply the $\Bbb R^+$--bundle whose smooth
sections are the metrics in the conformal class. Thus, the exact
Weyl--structures in conformal geometry correspond exactly to the
Levi--Civita connections of the metrics in the conformal class.

The reason for the names ``closed'' and ``exact'' becomes apparent,
once one studies the affine structures on the sets of closed and exact
Weyl--structures. So let us assume that $\si$ is a closed
Weyl--structure, and consider any other Weyl--structure $\hat\si$
corresponding to the section $\Up=(\Up_1,\dots,\Up_k)$ of $\Cal
A_1\oplus\dots\oplus\Cal A_k$. Now we can reinterpret theorem
\ref{3.12}(1) together with proposition \ref{3.9} as showing that the
set of Weyl--structures is an affine space over $\Om^1(M)$, in such a
way that fixing $\si$ the section $\Up$ corresponds to the one--form 
$\Up^{\si,\la}$ defined by 
$$
\Up^{\si,\la}(\xi)=\sum_{
\|\underline{j}\|+\ell=0} 
\frac{(-1)^{\underline{j}}}{\underline{j}!}\la'\left(\ad(\Up_k)^{j_k}\o
\dots\o\ad(\Up_1)^{j_1}(\xi_\ell)\right).
$$
This identification is obtained simply by pulling back the affine
structure on the space of principal connections on $\Cal L^\la$ to the
space of Weyl--structures. In particular, the change of the principal
connections $\tau$ and $\hat\tau$ on $\Cal L^\la$ induced by $\si$ and
$\hat\si$, respectively, is just given by
$\hat\tau=\tau+\Up^{\si,\la}$. But then their curvatures change simply
by $\hat\rho=\rho+d\Up^{\si,\la}$, so in particular if $\si$ is closed 
then $\hat\si$ is closed if and only if $d\Up^{\si,\la}=0$. Thus, in
the same way as Weyl--structures are affine over all one--forms,
closed Weyl--structures are affine over closed one--forms.

For exact Weyl--structures, the situation is even simpler. If $s$ and
$\hat s$ are two global sections of $\Cal L^\la$, then there is a
unique smooth function $f$ such that $\hat s(x)=e^{-f(x)}s(x)$. It is
then well known that the associated principal connections simply
change by $\hat\tau=\tau+df$, so exact Weyl--structures are affine
over the space of exact one--forms. 

\subsection{Remark}\label{3.14}
Another useful observation about exact Weyl geometries is related to the
identification of $\Cal L^\la$ with $\cg_0/\ker{\la}$ from \ref{3.8}(2). 
By the general
properties of classical $G$--structures, the sections of such bundles are in
bijective correspondence with reductions of the structure groups to
$\ker\la\subset G_0$. Thus the holonomy of the Weyl connections given by
closed Weyl structures is always at most $\ker\la$. In particular, in
$|1|$--graded cases the scaling element is unique up to scalar multiples, and 
the kernel of $\la$ is exactly the semisimple part of $G_0$. The same
observation is then true for the closed Weyl geometries locally.

\subsection{Normal Weyl--structures}\label{3.15}
Besides the rather obvious closed and exact Weyl--structures discussed above
there is a second kind of special Weyl--structures, the so--called {\em
normal Weyl--structures\/}.  In several respects, they are quite different
from closed and exact Weyl--structures. On one hand, they are ``more
canonical'' since their definition does not involve the choice of a bundle
of scales. On the other hand, in contrast to closed and exact
Weyl--structures, which always exist globally, normal Weyl--structures in
general exist only locally (over $M$). Their existence 
is closely related to the existence of normal coordinates
for parabolic geometries. This subject will be taken up elsewhere. We would
like to point out at this place that the existence of normal
Weyl--structures seems to be a new result even in the case of conformal
structures, where it significantly improves the result on the existence of
Graham normal coordinates, see \cite{Gra}. 

Since the Rho tensors give the information about
the difference of the covariant derivative with respect to 
the Weyl connection and the invariant derivative with respect to $\om$ along
the image of the chosen Weyl--structure $\si$, the ``normality'' we have in
mind will be described in terms of certain minimality of $\Rho$. More
explicitly, a Weyl--structure $\si$ will be called {\em normal at the point
$x\in M$} if it satisfies the properties imposed in Theorem \ref{3.16}.  
This Theorem also describes completely the freedom in the choice.

Recall that once we have chosen a Weyl--structure, we get an identification
of the tangent bundle with its associated graded vector bundle. Thus $TM$ is
associated to $\cg_0$ and so there is the induced linear Weyl connection on
$TM$. Since the Weyl--structure induces covariant derivatives on all
components of the associated graded of the tangent bundle, the Weyl
connection on $TM$ preserves the grading. For the same reason, we can form
covariant derivatives of the Rho--tensor, viewed as a one--form with values
in $T^*M\cong\Cal A_1\oplus\dots\oplus\Cal A_k$, which again preserve the
grading.

\begin{thm}\label{3.16}
Let $p:\Cal G\to M$ be a parabolic geometry with underlying
$G_0$--bundle $p_0:\Cal G_0\to M$ and let $\pi:\Cal G\to\Cal G_0$ be
the canonical projection. Let $x\in M$ be a point and let $u_0\in\Cal
G_0$ and $u\in\Cal G$ be points such that $\pi(u)=u_0$ and
$p_0(u_0)=x$. Then there exists an open neighborhood $U$ of $x$ in $M$
and a Weyl--structure $\si:p_0^{-1}(U)\to p^{-1}(U)$ such that
$\si(u_0)=u$ and the Rho--tensor $\Rho$ of $\si$ has the property that
for all $k\in\Bbb N$ the symmetrization over all $\xi_i$ of
$(\nabla_{\xi_k}\dots\nabla_{\xi_1}\Rho)(\xi_0)$ vanishes at $x$, so
in particular $\Rho(x)=0$. Moreover, this condition uniquely
determines the infinite jet of $\si$ in $u_0$.
\end{thm}
\begin{proof}
Consider the Cartan connection $\om$ on $\Cal G$. Since $\om$
restricts to a linear isomorphism, for each element $A\in\frak g$ we
get the constant vector field $\tilde A\in\frak X(\Cal G)$ defined by $\tilde
A(v)=\om(v)^{-1}(A)$, cf. \ref{2.4}. 
(Note that for $A\in\frak p$ this is just the
fundamental vector field.) In particular, we may consider the vector
fields $\tilde X$ for $X\in\frak g_-$. Now we can find a neighborhood
$V$ of zero in $\frak g_-$, such that for all $X\in V$ the flow of $X$
in the point $u$ exists up to time $t=1$. Define $\ph:V\to\Cal G$ by
$\ph(X)=\Fl^{\tilde X}_1(u)$. Since $T_up\o T_0\ph:\frak g_-\to T_xM$ is
obviously a linear isomorphism, we may assume (possibly shrinking $V$) 
that the maps $\ph$, $\pi\o\ph$ and $p\o\ph$ are all diffeomorphisms
onto their images, and we put $U=p(\ph(V))$. For a point $v_0\in
p_0^{-1}(U)$ there clearly exist unique elements $X\in V$ and $g\in
G_0$ such that $v_0=\pi(\ph(X))\cdot g$, and we define
$\si(v_0):=\ph(X)\cdot g$. Obviously, this defines a smooth
$G_0$--equivariant section $\si:p_0^{-1}(U)\to p^{-1}(U)$ and
$\si(u_0)=u$. 

Next, consider a tangent vector $\xi\in T_xM$, and its horizontal lift
$\xi^h\in T_{u_0}\Cal G_0$ with respect to the principal connection
$\si^*\om_0$. Since $\si^*\om_{\leq 0}$ defines a Cartan connection on
$p_0^{-1}(U)$ (see \ref{3.5}) we can extend $\xi^h$ uniquely to a
vector field $\tilde\xi^h$ such that $\si^*\om_{\leq 0}(\tilde\xi^h)$
is constantly equal to some $X\in\frak g_-$. Moreover, $\tilde\xi^h$
is projectable to a vector field $\tilde\xi$ on $U$ and it is exactly
the horizontal lift of $\tilde\xi$ (which also justifies the
notation). 

Now consider the flow line $c(t)=\Fl^{\tilde X}_t(u)=\ph(tX)$ in $\Cal
G$, which is defined for sufficiently small $t$. By construction, we
have $\si(\pi(c(t)))=c(t)$ for all $t$. But this implies that $T\si\cdot
(\pi\o c)'(t)=c'(t)$, so $\si^*(\om)((\pi\o c)'(t))$ is constantly
equal to $X$ and thus $(\pi\o c)(t)=\Fl^{\tilde\xi}_t(u_0)$. On the
other hand, from the construction it is clear that
$\om(c'(t))=X\in\frak g_-$, so if we consider the function $\Rho:\Cal
G_0\to L(\frak g_-,\frak p_+)$ describing the Rho--tensor, then
$\Rho(\pi(c(t)))(X)=0$ for all $t$. Consequently, all derivatives of
this curve in $t=0$ vanish. But since $\pi\o c$ is an integral curve
of $\tilde\xi^h$ these iterated derivatives exactly correspond to
iterated covariant derivatives of $\Rho$ in direction $\xi$ evaluated
at $\xi$. Thus, we obtain $(\nabla_\xi\dots\nabla_\xi\Rho)(\xi)=0$ for
any number of covariant derivatives. Using polarization, this implies
that the symmetrization of
$(\nabla_{\xi_k}\dots\nabla_{\xi_1}\Rho)(\xi_0)$ over all $\xi_i$
vanishes at $x$. 

To see that our condition fixes the infinite jet of the
Weyl--structure suppose that $\hat\si$ is another normal Weyl
structure with $\hat\si(u_0)=u$ and let $\Up=(\Up_1,\dots,\Up_k)$ be
the section of $\Cal A_1\oplus\dots\oplus\Cal A_k$ describing the change
from $\si$ to $\hat\si$. We want to show that the infinite jet of $\Up$
vanishes at $x=p(u)$. Since both Weyl--structures map $u_0$ to $u$,
we must have $\Up(x)=0$. Next, we know that
$\Rho(x)=\hat\Rho(x)=0$. Since all $\Up_i$ vanish in $x$, formula
\eqref{Rhotransf} from Proposition \ref{3.9} immediately shows that
this implies $\nabla\Up_i(x)=0$ for all $i=1,\dots,k$, so
$\nabla\Up(x)=0$. Now, $\hat\Rho=0$ and $\Up=0$. On one hand, it follows
that $\hat\nabla\Rho(x)=\nabla\Rho(x)$ and on the other hand that
$(\nabla_\eta\hat\Rho)(\xi)(x)=\nabla_\eta(\hat\Rho(\xi))(x)$. But
hitting formula \eqref{Rhotransf} from Proposition \ref{3.9} with
$\nabla_\eta$ and symmetrizing over $\xi$ and $\eta$, 
we always get terms involving some $\Up_i$ or
$\nabla_\eta\Up_i$ or $\nabla_\eta\Rho(\xi)$ which all vanish at $x$,
except for one term in the very last line, in which we get a
second covariant derivative of some $\Up_i$. So we see that the
symmetrizations of 
$\hat\nabla_\eta\hat\Rho(\xi)$ and $\nabla_\eta\nabla_\xi\Up$ coincide. Thus
vanishing of the symmetrization of the first covariant derivative
implies that the symmetrized second derivative of $\Up$ is zero, and
thus the two--jet of $\Up$ at $x$ must be zero. Iteratively, one
similarly sees that in the expression of an symmetrized iterated covariant
derivative of $\Rho$ we always get terms involving symmetrized 
iterated covariant
derivatives of $\Up_i$'s or $\Rho$'s except for one term coming from
the very last line of the transformation formula. As above, one then
concludes that vanishing of the symmetrization of the $k$--fold
covariant derivative of $\hat\Rho$ is equivalent to vanishing of the
symmetrization of the $(k+1)$--fold covariant derivative of $\Up$ and
thus to the $k+1$--jet of $\Up$ in $x$ being trivial.
\end{proof}

\section{Characterization of Weyl--structures}\label{4}
In the last section, we started with a Weyl--structure for a parabolic
geometry $(\cg\to M,\om)$ and we constructed several underlying objects on
the manifold $M$, see Figure 1 for an illustration. Now we are going to
characterize when general objects of that type actually come from a Weyl
structure. In the final stage, this will mean explicit conditions relating
the soldering form, linear connection and its torsion and curvature,
together with a procedure building the corresponding Rho--tensors. This 
is quite simple for irreducible geometries, where the
soldering form is fixed, and the whole condition prescribes uniquely the
torsion of a $G_0$--connection. The Rho--tensor is then given by a
simple formula in terms of the curvature, see Example \ref{4.8}
below. Of course, the same story gets much more complicated for the general
$|k|$--graded case. The main step is done in Theorem \ref{4.4} and then a
detailed analysis of the curvature fulfills our goal.

Throughout this section we restrict to the case of regular parabolic geometries
associated to a $|k|$--graded semisimple Lie algebra $\frak g$ such
that $H^1(\frak g_-,\frak g)$ is concentrated in homogeneous degrees
$\leq 0$, i.e.~such that none of the simple $|k_i|$--graded ideals is
of one of the two types mentioned  in \ref{2.8}. In the case that such
ideals are present, a similar characterization is possible, but the
conditions are more complicated to formulate. 

\subsection{Definition}\label{4.1}
Let $p_0:\Cal G_0\to M$ be a regular infinitesimal flag structure, see
\ref{2.7}. A {\em Weyl--form\/} for $M$ is a one--form
$\tau\in\Om^1(\Cal G_0,\frak g)$ which is $G_0$--equivariant,
i.e.~$(r^g)^*\tau=\Ad(g^{-1})\o\tau$ for all $g\in G_0$, reproduces
the generators of fundamental vector fields, i.e.~$\tau(\ze_A)=A$ for
all $A\in\frak g_0$ and has the property that for each $i<0$ the
restriction of $\tau_i$ to $T^i\Cal G_0$ coincides with the $\frak
g_i$--component of the frame form of degree one on $\Cal G_0$ induced
by the regular infinitesimal flag structure, see \ref{2.7} and \ref{3.5}.

\medskip

By \ref{3.3} and \ref{3.4}, 
for any Weyl--structure $\si:\Cal G_0\to\Cal G$, the
pullback  $\si^*\om$ is a Weyl--form for $M$. As
in \ref{3.4}, the
condition of the restriction of $\tau_i$ to $T^i\Cal G_0$, $i<0$, 
means on $M$
exactly that the restriction of $\tau_i$ to $T^iM$ coincides with the
canonical projection $T^iM\to\Cal A_i$. In particular, this implies
that $\tau_-=\tau_{-k}\oplus\dots\oplus\tau_{-1}$ induces a linear
isomorphism $T_u\Cal G_0/V_u\Cal G_0\cong\frak g_-$, and thus
$\tau_{\leq 0}$ is a Cartan connection on $\Cal G_0$. Completely
parallel to the development in \ref{3.3}--\ref{3.5}
we can equivalently interpret
a Weyl--form for $M$ as a one form $\tau_-\in\Om^1(M,\Cal
A_{-k}\oplus\dots\oplus\Cal A_{-1})$ inducing an isomorphism between
$TM$ and its associated graded bundle, plus a principal connection
$\tau_0\in\Om^1(\Cal G_0,\frak g_0)$ on $\Cal G_0$, plus a Rho--tensor
$\Rho=\Rho^\tau\in\Om^1(M,\Cal A_1\oplus\dots\oplus\Cal A_k)$, so a
Weyl--form essentially consists of objects living on $M$.

\subsection{Weyl--curvature}\label{4.2}
Next, for a Weyl--form $\tau$ for $M$, we define the {\em
Weyl--curvature\/} $W$ of $\tau$. As a $\frak g$--valued two form on $\Cal G_0$, it is defined by
$$
W(\xi,\eta)=d\tau(\xi,\eta)+[\tau(\xi),\tau(\eta)].
$$
From the fact that $\tau$ is $G_0$--equivariant and reproduces the
generators of fundamental vector fields, one immediately concludes that
$W$ is horizontal and $G_0$--equivariant, so it descends to an $\Cal
A$--valued two form on $M$. Taking into account the identification of
$TM$ with $\Cal A_-$, we can also view $W$ as a section of
$L(\La^2\Cal A_-,\Cal A)$.

Finally note that any section $\Ph$ of $L(\La^2\Cal A_-,\Cal A)$ can be
split according to homogeneous degrees. We denote by $\Ph^{(\ell)}$ the
homogeneous part of degree $\ell$, i.e.~$\Ph^{(\ell)}(\xi,\eta)\in
\Cal A_{i+j+\ell}$ for sections $\xi$ of $\Cal A_i$ and $\eta$ of $\Cal
A_j$ with $i,j<0$. 

\begin{lem*}
Let $p_0:\Cal G_0\to M$ be a regular infinitesimal flag structure.
Then any Weyl--form
$\tau\in\Om^1(\Cal G_0,\frak g)$ has the property that $W^{(\ell)}=0$
for all $\ell\leq 0$.  
\end{lem*}
\begin{proof}
Consider $\xi\in\Ga(\Cal A_i)$ and $\eta\in\Ga(\Cal A_j)$, for
$i,j<0$. Then $\tau_n(\xi)=0$ for $n<i$ and $\tau_m(\eta)=0$ for
$m<j$, so for $\ell<0$ if $m+n=i+j+\ell$ then
$[\tau_n(\xi),\tau_m(\eta)]=0$. Thus, in this case, the definition of
$W^{(\ell)}(\xi,\eta)$ can be rewritten as
$W^{(\ell)}(\xi,\eta)=d\tau_{i+j+\ell}(\xi,\eta)=
-\tau_{i+j+\ell}([\xi,\eta])$. By definition of a Weyl--form,
$W^{(\ell)}(\xi,\eta)$ thus equals the class of the bracket
$-[\xi,\eta]$ in $T^{i+j+\ell}M/T^{i+j+\ell+1}M$. But according to 
\ref{2.7}, we in particular know
that the bracket of any section of $T^iM$ with a section of $T^jM$
lies in $T^{i+j}M$, so since $\ell<0$, we must have $W^{(\ell)}=0$. 

Next, for $\ell=0$, we can write 
$$
W^{(0)}(\xi,\eta)=d\tau_{i+j}(\xi,\eta)+\{\xi,\eta\}=
-\tau_{i+j}([\xi,\eta])+\{\xi,\eta\}.
$$
Again, $\tau_{i+j}([\xi,\eta])$ is just the class of the bracket in
$T^{i+j}M/T^{i+j+1}M$ and so the vanishing of $W^{(0)}$ is just the
remaining part of the definition of regular infinitesimal flag structures,
see \ref{2.7}.
\end{proof}

\subsection{Definition}\label{4.3}
Let $p_0:\Cal G_0\to M$ be a regular infinitesimal flag structure.
Then a Weyl--form
$\tau\in\Om^1(\Cal G_0,\frak g)$ is called {\em normal\/} if and only
if its Weyl--curvature $W\in\Ga(L(\La^2\Cal A_-,\Cal A))$ satisfies
$\partial^*(W)=0$, where $\partial^*:L(\La^2\Cal A_-,\Cal A))\to
L(\Cal A_-,\Cal A)$ is the bundle map induced by the Lie algebra
codifferential, see \ref{2.13}.

\begin{thm}\label{4.4}
Let $(p:\Cal G\to M,\om)$ be a regular parabolic geometry and let
$p_0:\Cal G_0\to M$ be the underlying regular infinitesimal flag structure. 
Then a Weyl--form $\tau\in\Om^1(\Cal G_0,\frak g)$ for $M$ is coming
from some Weyl--structure $\si:\Cal G_0\to\Cal
G$, i.e. $\tau=\si^*\om$,  if and only if $\tau$ is normal.
\end{thm}
\begin{proof}
First we show that for any Weyl--structure $\si:\Cal G_0\to\Cal G$ the
Weyl--form $\si^*\om$ is normal. By  the definition in \ref{4.2} the
Weyl--curvature $W$ is a $\frak g$--valued two--form on $\Cal G_0$, 
given by 
$$
W(\xi,\eta)=d\si^*\om(\xi,\eta)+[\si^*\om(\xi),\si^*\om(\eta)]=
\si^*(d\om+\tfrac{1}{2}[\om,\om])(\xi,\eta). 
$$
Thus, $W$ is simply the pullback along $\si$ of the curvature of the
Cartan connection $\om$ on $\Cal G_0$. By definition of a normal
parabolic geometry, this curvature is $\partial^*$--closed, so the same
is true for $W$.

Now, let us assume that we have given an arbitrary normal Weyl--form
$\tau\in\Om^1(\Cal G_0,\frak g)$. Moreover, let us choose any bundle
$\Cal L^\la$ of scales for the parabolic geometry in question. Since
$\tau_0$ is a principal connection on $\Cal G_0$, it induces a
principal connection on $\Cal L^\la$, which by Theorem \ref{3.12} in
turn gives rise to a unique Weyl--structure $\si:\Cal G_0\to\Cal G$
such that the connection on $\Cal L^\la$ induced by $\si^*\om_0$
coincides with the connection induced by $\tau_0$. We claim that
$\tau=\si^*\om$, which will conclude the proof. 

Consider the difference $\tau-\si^*\om\in\Om^1(\Cal G_0,\frak g)$. For
$i<0$, we know from our assumptions that both $\tau_i$ and
$\si^*\om_i$ coincide on $T^i\Cal G_0$ with the frame form of degree
one. In particular, the difference $\tau_i-\si^*\om_i$ vanishes on
$T^i\Cal G_0$ for all $i<0$. 
Since $T^0\Cal G_0$ is just the vertical bundle of
$\Cal G_0$ and since both $\tau_0$ and $\si^*\om_0$ are principal
connections on $\Cal G_0$, we see that $\tau_0-\si^*\om_0$ vanishes on
$T^0\Cal G_0$, too. Finally, if we put $T^i\Cal G_0$ to be the zero section
for $i>0$, then $\tau_i-\si^*\om_i$ vanishes on $T^i\Cal G_0$ for all
$i=-k,\dots, k$. Let us inductively assume that $\tau_i-\si^*\om_i$
vanishes on $T^{i-n+1}\Cal G_0$ for all $i$ and some $n$. 

Then consider the restriction of $\tau_i-\si^*\om_i$ to 
$T^{i-n}\Cal G_0$, which can be viewed as a map $T^{i-n}\Cal
G_0/T^{i-n+1}\Cal G_0\to\frak g_i$. For each $i$ such that $i-n\leq
0$, the forms $\tau_{i-n}$ and $\si^*\om_{i-n}$ coincide on
$T^{i-n}\Cal G_0$ and induce an isomorphism$T^{i-n}\Cal
G_0/T^{i-n+1}\Cal G_0\to \Cal G_0\x\frak g_{i-n}$. Consequently, we
get a unique map $\Ph:\Cal G_0\to L(\frak g_-,\frak g)$ which has
values in maps homogeneous of degree $n$, such that
$(\tau_i-\si^*\om_i)(\xi)=\Ph(\tau_{i-n}(\xi))$ for all $\xi\in
T^{i-n}\Cal G_0$. 

Next, let $W^{(n)}$ be the homogeneous component of degree $n$
of the Weyl--curvature of $\tau$ viewed as a function $\Cal G_0\to
L(\La^2\frak g_-,\frak g)$ (having values in the maps homogeneous of
degree $n$), and let $\tilde W^{(n)}$ be the corresponding object for
$\si^*\om$. We claim that for all $X$, $Y\in \frak g_-$
\begin{align}\label{W-relation}
\tilde W^{(n)}(X,Y)&=W^{(n)}(X,Y)-[X,\Ph(Y)]+[Y,\Ph(X)]+\Ph([X,Y])=\\
&=W^{(n)}(X,Y)-(\partial\o\Ph)(X,Y). 
\notag\end{align}

Let us postpone the proof of this claim and 
assume it is true for a while. 
Since both $W^{(n)}$ and $\tilde W^{(n)}$ are $\partial^*$--closed,
this implies $\partial^*\o\partial\o\Ph=0$, which implies
$\partial\o\Ph=0$ since $\partial$ and $\partial^*$ are adjoint, see
\ref{2.5}. Since $H^1(\frak g_-,\frak g)$ is concentrated in
non-positive degrees of homogeneity, this implies $\Ph=0$ for $n>k$ and
$\Ph(X)=[Z,X]$ for some smooth $Z:\Cal G_0\to\frak g_n$ for $n\leq
k$. But in the latter case, the proof of Theorem \ref{3.12}(1) shows that
since $\tau_0$ and $\si^*\om_0$ induce the same principal connection on
$\Cal L^\la$, we must have $Z=0$, and thus $\Ph=0$. Hence, $\tau_i$ and
$\si^*\om_i$ coincide on $T^{i-n}\Cal G_0$ for all $i<n$, for $i=n$ this
follows since $n>0$ and thus both
$\tau_n$ and $\si^*\om_n$ are horizontal, while for $i>n$ it is
trivially satisfied. Thus the result follows by induction.  

So we are left with the proof of \eqref{W-relation} only. Let us fix
$X\in\frak g_i$, $Y\in\frak g_j$, $i,j<0$.
By definition, 
$$
W^{(n)}(u)(X,Y)=d\tau_{i+j+n}(\tau_{\leq 0}^{-1}(X),\tau_{\leq
0}^{-1}(Y))+[\tau(\tau_{\leq 0}^{-1}(X)),\tau(\tau_{\leq
0}^{-1}(Y))]_{i+j+n}, 
$$
where the index
in the bracket means that we just have to take the component in $\frak
g_{i+j+n}$. For $\tilde W^{(n)}$ we get the analogous formula with all
$\tau$'s replaced by $\si^*\om$. 

Next, observe that both $\tau_{\leq 0}^{-1}(X)$ and $\si^*\om_{\leq
0}^{-1}(X)$ lie in $T^i\Cal G_0\subset T^{i+j}\Cal G_0$ and similarly
for $Y$. From above, we know that
$\si^*\om_{i+j+n}(\xi)=\tau_{i+j+n}(\xi)-\Ph(\tau_{i+j}(\xi))$ for all
$\xi$ in $T^{i+j}\Cal G_0$. Taking the exterior derivative of this
equation and keeping in mind that $\tau_{i+j}$ vanishes on $T^i\Cal
G_0$ and $T^j\Cal G_0$, we see that for $\xi\in T^i\Cal G_0$ and
$\eta\in T^j\Cal G_0$ we get 
$$
d\si^*\om_{i+j+n}(\xi,\eta)=d\tau_{i+j+n}(\xi,\eta)-
\Ph(d\tau_{i+j}(\xi,\eta)). 
$$
Since $W^{(0)}=0$, the second term (including the $-$ sign) can be
rewritten as $\Ph([\tau_i(\xi),\tau_j(\eta)])$, and we may as well
replace $\tau$ by $\si^*\om$ in this expression. Thus, we see that
\begin{align*}
\tilde W^{(n)}(X,Y)&=d\tau_{i+j+n}(\si^*\om_{\leq 0}^{-1}(X),\si^*\om_{\leq
0}^{-1}(Y))+\Ph([X,Y])+\\
&+[\si^*\om(\si^*\om_{\leq 0}^{-1}(X)),
\si^*\om(\si^*\om_{\leq 0}^{-1}(Y))]_{i+j+n}.
\end{align*}
Now we have to distinguish a few cases: Let us first assume that
$i+n>0$. Then $\si^*\om_{\leq 0}^{-1}(X)=\tau_{\leq 0}^{-1}(X)$, and
$\si^*\om(\si^*\om_{\leq 0}^{-1}(X))=\tau(\tau_{\leq
0}^{-1}(X))-\Ph(X)$, and $\Ph(X)\in\frak g_{i+n}\subset\frak
p_+$. In particular, this implies that 
\begin{multline*}
[\si^*\om(\si^*\om_{\leq 0}^{-1}(X)),
\si^*\om(\si^*\om_{\leq 0}^{-1}(Y))]_{i+j+n}=\\
=[\tau(\tau_{\leq 0}^{-1}(X)),\si^*\om(\si^*\om_{\leq
0}^{-1}(Y))]_{i+j+n}-[\Ph(X),Y].
\end{multline*}
Secondly, if $i+n=0$ then $\Ph(X)\in\frak g_0$, and thus
$\si^*\om_{\leq 0}^{-1}(X)=\tau_{\leq 0}^{-1}(X)+\ze_{\Ph(X)}$. The
infinitesimal version of equivariance of $\tau_{i+j+n}$ then implies
that 
$$
d\tau_{i+j+n}(\si^*\om_{\leq 0}^{-1}(X),\si^*\om_{\leq 0}^{-1}(Y))=
d\tau_{i+j+n}(\tau_{\leq 0}^{-1}(X),\si^*\om_{\leq 0}^{-1}(Y))-[\Ph(X),Y],
$$
since $i+j+n=j$ in this case. On the other hand both
$\si^*\om(\si^*\om_{\leq 0}^{-1}(X))$ and $\tau(\tau_{\leq
0}^{-1}(X))$ in this case are congruent to $X$ modulo $\frak p_+$, so 
$$
[\si^*\om(\si^*\om_{\leq 0}^{-1}(X)),
\si^*\om(\si^*\om_{\leq 0}^{-1}(Y))]_{i+j+n}=
[\tau(\tau_{\leq 0}^{-1}(X)),
\si^*\om(\si^*\om_{\leq 0}^{-1}(Y))]_{i+j+n}.
$$
Finally, suppose that $i+n<0$, so $\Ph(X)\in\frak g_{i+n}\subset\frak
g_-$. Then $\si^*\om_{\leq 0}^{-1}(X)$ is congruent to $\tau_{\leq
0}^{-1}(X+\Ph(X))$ modulo $T^{i+n+1}\Cal G_0$. Since the bracket of
a section of this subbundle with a section of $T^j\Cal G_0$ is a
section of $T^{i+j+n+1}\Cal G_0$ and $\tau_{i+j+n}$ vanishes on the
latter subbundle, we conclude that 
\begin{align*}
d\tau_{i+j+n}(\si^*\om_{\leq 0}^{-1}(X),\si^*\om_{\leq 0}^{-1}(Y))&=
d\tau_{i+j+n}(\tau_{\leq 0}^{-1}(X),\si^*\om_{\leq 0}^{-1}(Y))+\\
&+d\tau_{i+j+n}(\tau_{\leq 0}^{-1}(\Ph(X)),\si^*\om_{\leq 0}^{-1}(Y)).
\end{align*}
Since $W^{(0)}=0$, the last term can be rewritten as $-[\Ph(X),Y]$. As
above, both $\si^*\om(\si^*\om_{\leq 0}^{-1}(X))$ and $\tau(\tau_{\leq
0}^{-1}(X))$ are congruent to $X$ modulo $\frak p_+$, so again the
bracket term makes no problem.

Hence we see, that in any case we get 
\begin{align*}
\tilde W^{(n)}(X,Y)&=d\tau_{i+j+n}(\tau_{\leq 0}^{-1}(X),\si^*\om_{\leq
0}^{-1}(Y))+[\tau(\tau_{\leq 0}^{-1}(X)),
\si^*\om(\si^*\om_{\leq 0}^{-1}(Y))]_{i+j+n}-\\
&-[\Ph(X),Y]+\Ph([X,Y]). 
\end{align*}
Doing the same changes to $Y$ instead of $X$ we obtain the required equality
\eqref{W-relation}, and the whole proof of the theorem is finished.
\end{proof}

\subsection{Remark}\label{4.5}
If one does not assume that $H^1(\frak g_-,\frak g)$ is concentrated
in non-positive degrees, i.e.~if one allows $\frak g$ to contain one of
the two simple factors mentioned in \ref{2.8}, then $H^1(\frak
g_-,\frak g)$ is concentrated in homogeneous degrees less or equal to
one. Thus, the above proof shows that $\tau=\si^*\om$ if $\tau$ is
normal and has the property that the restrictions of $\tau_i$ and
$\si^*\om_i$ to $T^{i-1}\Cal G_0$ coincide for all $i$. This condition
is then fairly simple to interpret for any concrete choice of such
structure. 

\subsection{}\label{4.6}
In the proof of Theorem \ref{4.4}, we observed that for a
Weyl--structure $\si:\Cal G_0\to\Cal G$ the Weyl--curvature $W$ of the
Weyl--form $\si^*\om$ is exactly the pullback along $\si$ of the
curvature $\ka$ of the normal Cartan connection $\om$ on $\Cal
G$. This allows us to compute the change of the Weyl--curvature under
a change of the Weyl--structure. Suppose that $\hat\si$ is another
Weyl--structure and $\Up=(\Up_1,\dots,\Up_k)$ is the smooth section of
$\Cal A_1\oplus\dots\oplus\Cal A_k$ describing the change from $\si$
to $\hat\si$, see Proposition \ref{3.2}, i.e.
$$
\hat\si(u)=\si(u)\exp(\Up_1(u))\dots\exp(\Up_k(u)).
$$
Equivariance of the Cartan connection $\om$ immediately implies that
the curvature $\ka$ is equivariant, i.e.~viewing $\ka$ as a two form
on $\Cal G$ with values in $\frak g$, we have $\ka(v\cdot
g)(Tr^g\cdot\xi,Tr^g\cdot\eta)=\Ad(g^{-1})(\ka(v)(\xi,\eta))$ for
$g\in P$ and $\xi,\eta\in T_v\Cal G$. Putting $v=\si(u)$ and
$g=\exp(\Up_1(u))\dots\exp(\Up_k(u))$, we see from the proof of
Proposition \ref{3.9} that for $\xi\in T_u\Cal G_0$ the element
$T_u\hat\si\cdot\xi$ is congruent to $Tr^gT_u\si\cdot\xi$ modulo
vertical elements, which are killed by the curvature anyhow. Thus,
viewing $W$ and $\hat W$ as $\frak g$--valued two forms on $\Cal G_0$,
we get $\hat W(\xi,\eta)=\Ad(g^{-1})(W(\xi,\eta))$. Moreover, to get
the interpretation of our two Weyl curvatures $W$ and $\hat W$ as $\Cal A$--valued two forms
on $M$, we just have to apply the above definition to lifts of vector
fields on $M$, and the result is independent of the choice of the
lifts since $W$ is horizontal. Keeping in mind that the Lie--bracket
in $\frak g$ corresponds to the algebraic bracket of sections of $\Cal
A$ and expanding the exponentials in $\Ad(g^{-1})$ as in the proof of
Proposition \ref{3.9} we arrive (with notation as in \ref{3.9}) at
\begin{equation}\label{Weyltransf}
\hat W_i(\xi,\eta)=\sum_{
\|\underline{j}\|+\ell=i}
\frac{(-1)^{\underline{j}}}{\underline{j}!}\ad(\Up_k)^{j_k}\o\dots
\o\ad(\Up_1)^{j_1}(W_\ell(\xi,\eta)).
\end{equation}
From this formula, one can also derive a formula describing the change
of $W$ viewed as a section of $L(\La^2\Cal A_-,\Cal A)$ taking into
account the change of the identification of $TM$ with $\Cal A_-$
described by \eqref{splittransf} in Proposition \ref{3.9}, 
and thus a formula for the change
of the individual homogeneous components $W^{(\ell)}$. The only point
that is important for us here is that {\em the homogeneous component
$W^{(1)}$ of degree one is actually independent of $\si$}. This can be
immediately verified from the above formula, taking into account that
$W^{(\ell)}=0$ for all $\ell\leq 0$. 


\subsection{Remark}\label{4.7}
The results obtained so far in principle allow to give a description
of the Cartan bundle and the Cartan connection completely in terms of
data on the manifold $M$. More precisely, if we start from a regular
infinitesimal flag structure underlying
some parabolic geometry, then we may proceed as follows: Choose a
scaling element $E_\la\in\frak z(\frak g_0)$, and consider the
corresponding homomorphism $\la:G_0\to\Bbb R^+$ described in the proof
of Proposition \ref{3.8}. Then form $\Cal L_\la=\Cal G_0\x_{G_0}\Bbb
R^+$. From Theorem \ref{3.12}(2) we then know that the Cartan bundle
$\Cal G$ is just the pullback of the bundle of principal connections
on $\Cal L_\la$, and we have a description of the principal
action. Moreover, a choice of a principal connection on $\Cal L^\la$
is just the choice of a global section of the bundle of connections,
so its pullback is a smooth $G_0$--equivariant section $\si:\Cal
G_0\to\Cal G$. Any Cartan connection $\om$ on $\Cal G$ is uniquely
determined by its pullback $\si^*\om$ by equivariance. Thus,
describing the canonical normal Cartan connection on $\Cal G$ is
equivalent to finding a normal Weyl--form on $\Cal G_0$ which induces
a given connection on $\Cal L^\la$.

\subsection{Example}\label{4.8}
Let us look more closely at the irreducible parabolic geometries. 
Here the regular infinitesimal flag structures are just $G_0$--structures 
on $M$ in the sense of classical G--structures. 
The Weyl forms are $\tau=\tau_{-1}+\tau_0+\tau_1$ where 
$\tau_{-1}:T\cg_0\to \frak g_{-1}$ is the fixed soldering form for $M$, 
$\tau_0$ is any linear connection on $M$ belonging to the fixed 
$G_0$--structure and $\tau_1$ is any one--form in $\Om^1(M;T^*M)$. Now, 
$$
W_{-1} = d\tau_{-1} + [\tau_{-1},\tau_0]
,$$
i.e. the torsion of the connection $\tau_0$. The individual components of $W$ have homogeneities one,
two, and three and so they have to be $\partial^*$--closed separately. The condition $\partial^*W_{-1}=0$ means that the torsion of $\tau_0$ is harmonic and this is the part of $W$ independent of the choice of the Weyl--structure. Next,
$$
W_0=d\tau_0+\tfrac12[\tau_0,\tau_0] + [\tau_{-1},\tau_1]
$$
which is the curvature $R$ of the connection $\tau_0$ plus some additional
term. The co--closedness of $W_0$ imposes a condition on the choice of
$\tau_1$, while $\partial^*W_1$ always vanishes since its values are in the
trivial vector space.

We shall see later that the resulting system of equations for the tensor
$\tau_1$ is always solvable, except for the projective structures (where the
first cohomology is concentrated in degree one). Moreover, we shall prove an
explicit algebraic formula for the necessary choice for the Rho--tensor:
$\tau_1=\square^{-1}\partial^*R$. Expanding this formula in the case of the
conformal (pseudo) Riemannian geometry, we obtain the well known Rho--tensor
used heavily by many authors since the beginning of this century, while
$d\tau_1$ happens to be exactly another well known tensor, the Cotton--York
tensor.

As mentioned above, this computation may be understood as an alternative for
the explicit construction of the canonical Cartan connection for all
irreducible parabolic geometries. 

\subsection{Total curvature}\label{4.9}
The explicit construction of a normal Weyl--form depends a lot on
the structure in question, a detailed treatment in the case of
partially integrable almost CR--structures of hypersurface type will
appear in \cite{Cap}. Here we just describe the basic ingredient of
this procedure. The upshot of this is that the condition on a
Weyl--form $\tau$ being normal can be step by step reduced to a
condition on $\tau_{\leq 0}$ only, at the same time computing step by
step the components of the Rho--tensor $\Rho=\tau_+$.

The first step in this direction is to 
replace the Weyl curvature of a Weyl--form $\tau$ by 
2--forms defined by splitting the structure equations for $\tau$.  
The curvature of the Cartan connection $\tau_{\le0}$
is the 2--form  $K_{\leq 0}\in\Om^2(\cg_0,\frak g_-)$ given by
$$
K_{\leq 0}(\xi,\eta)=d\tau_{\leq 0}(\xi,\eta)+[\tau_{\leq
0}(\xi),\tau_{\leq 0}(\eta)].
$$
On the other hand, we define the 2--form $K_+\in\Om^2(\cg_0,\frak p_+)$ by 
$$
K_+(\xi,\eta)=d\tau_+(\xi,\eta)+[\tau_+(\xi),\tau_+(\eta)].
$$
Motivated by conformal geometry, we call $K_+$ the {\em
Cotton--York--tensor} associated to the Weyl--form $\tau$. We write
$K=K_{\le0}+K_+$ and we call it the {\em (total) curvature of $\tau$}. 
Since
$\tau_{\leq 0}$ is a Cartan connection, it is well known that its curvature
is horizontal and $G_0$--equivariant, so it can be viewed as a two form on
$M$, with values in the bundle $\Cal A_{-k}\oplus\dots\oplus\Cal A_0$. On
the other hand, since $\tau_+$ is by assumption $G_0$--equivariant and
horizontal, the part $K_+$ descends to $M$, too. Finally, taking into
account the isomorphism
$TM\cong\Cal A_-=\Cal A_{-k}\oplus\dots\oplus\Cal A_{-1}$, we can finally
view $K$ as a smooth section of the bundle $L(\La^2\Cal A_-,\Cal A)$ over
$M$.

The reason for introducing this curvature is that it is more closely related
to usual invariants of the Weyl--form than the Weyl--curvature, cf. Example
\ref{4.8}. On the other hand, we shall see that there still is a simple
relation between curvature and Weyl--curvature.

To get explicit expressions for the components of $K$, recall that the
component $\tau_0$ of any Weyl--form $\tau$ is a principal
connection on $\Cal G_0$, and thus induces a linear connection
$\nabla$ on each of the bundles $\Cal A_i$. Let us also recall that $\tau_i$
are identified with forms $\Om^1(M;\ca_i)$ for all negative $i$.

\begin{prop}\label{4.10}
Let $p_0:\Cal G_0\to M$ be a regular infinitesimal flag structure, let
$\tau\in\Om^1(\Cal G_0,\frak g)$ be a Weyl--form for $M$, and let $K$
be its total curvature, viewed as an $\Cal A$--valued two form on $M$ with
$\Cal A_\ell$--component $K_\ell$. Then for all vector fields $\xi$ and
$\eta$ on $M$ we have:\newline
(1) $K_\ell(\xi,\eta)=\nabla_\xi(\tau_\ell(\eta))-
\nabla_\eta(\tau_\ell(\xi))-\tau_\ell([\xi,\eta])+
\sum_{\substack{i,j<0\\i+j=\ell}}\{\tau_i(\xi),\tau_j(\eta)\}$, for
$\ell<0$.\newline
(2) For $\zeta\in\Cal A_m$ we get
$\{K_0(\xi,\eta),\zeta\}=R_m(\xi,\eta)(\zeta)$, where $R_m$ is the
curvature of the linear connection $\nabla$ on $\Cal A_m$.

Moreover, if we view $K$ as a section of $L(\La^2\Cal A_-,\Cal A)$ and
consider $\ell>0$, then the homogeneous component $K^{(\ell)}$ of $K$ depends
only on the restrictions of $\tau_i$ to $T^{i-\ell}\Cal G_0$ for all
$i\leq 0$ and on the restrictions of $\tau_i$ to $T^{i-\ell+1}\Cal
G_0$ for $i>0$.
\end{prop}
\begin{proof}
By definition, for $\ell<0$ the function $\Cal G_0\to\frak g_\ell$
corresponding to $K_\ell(\xi,\eta)$ is given by 
\begin{multline*}
d\tau_\ell(\xi^h,\eta^h)+\sum_{i,j\leq 0,i+j=\ell}
[\tau_i(\xi^h),\tau^j(\eta^h)]=\\
=\xi^h\cdot\tau_\ell(\eta^h)-\eta^h\cdot\tau_\ell(\xi^h)-
\tau_\ell([\xi^h,\eta^h])+\sum_{i,j\leq 0,i+j=\ell}
[\tau_i(\xi^h),\tau_j(\eta^h)],
\end{multline*}
where the superscript $h$ denotes the horizontal lift with respect to
the principal connection $\tau_0$. But now $\tau_\ell(\eta^h):\Cal
G_0\to\frak g_\ell$ is exactly the smooth function corresponding to
the section $\tau_\ell(\eta)$ of $\Cal A_\ell$, so the function
$\xi^h\cdot\tau_\ell(\eta^h)$ corresponds to
$\nabla_\xi(\tau_\ell(\eta))$ and similarly for the second term. On
the other hand, $[\xi^h,\eta^h]$ is a lift of the vector field
$[\xi,\eta]$, so since $\tau_\ell$ is horizontal for $\ell<0$, we see
that the function $\tau_\ell([\xi^h,\eta^h])$ corresponds to the
section $\tau_\ell([\xi,\eta])$ of $\Cal A_\ell$. Finally, for the
last sum one only has to take into account that $\tau_0$ vanishes on
horizontal lifts and the bracket in $\frak g$ corresponds to the
algebraic bracket on $\Cal A$.

If $\ell=0$, the definition of $K_0$ reduces to
$d\tau_0(\xi^h,\eta^h)$ and this exactly represents the curvature of
the principal connection $\tau_0$, so the result follows immediately,
taking into account that the action of $\frak g_0$ on $\frak g_m$
is given by the Lie bracket in $\frak g$ and thus corresponds to the
algebraic bracket $\Cal A_0\x\Cal A_m\to\Cal A_m$.

To verify the statements about homogeneous degrees, take sections
$\xi$ of $\Cal A_i$ and $\eta$ of $\Cal A_j$, and let $\tilde\xi$ be
the (unique) section of $T^iM$ such that $\tau_n(\tilde\xi)=0$ for all
$i<n<0$, $\tau_i(\tilde\xi)=\xi$, 
and similarly for $\tilde\eta$. Then for $\ell>0$,
$K^{(\ell)}(\xi,\eta)=K_{i+j+\ell}(\tilde\xi,\tilde\eta)$. If
$i+j+\ell<0$, then the above formula just gives us
$$
\delta^0_{i+\ell}\nabla_{\tilde\xi}\eta-
\delta^0_{j+\ell}\nabla_{\tilde\eta}\xi
-\tau_{i+j+\ell}([\tilde\xi,\tilde\eta]).
$$
This is completely independent of the components $\tau_n$ for
$n>0$. If we allow a change of $\tau$ without changing the restriction
of $\tau_n$ to $T^{n-\ell}$ for all $n\leq 0$, then this means that
$\tilde\xi$ is changed at most by a section of $T^{i+\ell+1}M$. In
particular, if the first term in the above expression actually occurs,
i.e.~$i+\ell=0$ then $\tilde\xi$ is fixed, and moreover, since the
restriction of $\tau_0$ to $T^{-\ell}\Cal G_0=T^i\Cal G_0$ is fixed,
also the covariant derivative is fixed. Similarly one analyzes the
second term. Finally, the last term depends only on the restriction of
$\tau$ since the bracket of a section of $T^{i+\ell+1}M$ with a
section of $T^jM$ is a section of $T^{i+j+\ell+1}M$ and this subbundle
lies in the kernel of $\tau_{i+j+\ell}$. 

If $i+j+\ell=0$, then
$K^{(\ell)}(\xi,\eta)=d\tau_0((\tilde\xi)^h,(\tilde\eta)^h)$, and as
above, we see that $(\tilde\xi)^h$ and $(\tilde\eta)^h$ depend only on
the appropriate restriction of $\tau$. Moreover, the bracket
$[(\tilde\xi)^h,(\tilde\eta)^h]$ by construction is a section of
$T^{i+j}\Cal G_0$, so the whole expression depends only on the
restriction of $\tau_0$ to $T^{i+j}\Cal G_0=T^{0-\ell}\Cal G_0$. 

Finally, we have to consider the case $i+j+\ell>0$, so we are dealing
with a component of $K$ having values in $\Cal A_+$. As before, one
verifies that all extensions and horizontal lifts depend only on the
appropriate restrictions of $\tau_{\leq 0}$, so what remains to be
discussed is the dependence on $\Rho$. But viewing $\Rho$ as a section
of $L(\Cal A_-,\Cal A_+)$, the statement to be proved reduces to the
fact that a homogeneous component of $K$ depends only on homogeneous
components of $\Rho$ of strictly smaller degree. But this is obvious
from the definition of $K_+$.
\end{proof}

\subsection{Remark}\label{4.11}
The previous Proposition reveals that the $\ca_-$--components of the total
curvature give exactly the torsion of the linear connection $\tau_0$
corrected by the algebraic contribution of the Lie bracket in $\frak g_-$,
while the component $K_0$ is just the standard curvature of $\tau_0$. For
a normal Weyl form $\tau$ this means (using Proposition \ref{4.12} below)
that the torsion of $\tau_0$ has the
algebraic bracket as its homogeneous component of degree zero, no components
of negative degrees, and some positive degree components. The torsion 
component of degree one is an invariant of the parabolic structure in
question.

The key point in the further analysis is that while the total curvature of a
Weyl--form is much easier to relate to the underlying structure than
its Weyl--curvature, there is the quite simple relation between
them described in the next Proposition.

\begin{prop}\label{4.12}
Let $\tau\in\Om^1(\Cal G_0,\frak g)$ be a Weyl--form for $M$, let
$\Rho\in\Ga(L(\Cal A_-,\Cal A_+))$ be its Rho--tensor, and let $K,W\in
\Ga(L(\La^2\Cal A_-,\Cal A))$ be its total curvature and its
Weyl--curvature, respectively. Then 
$$
W(\xi,\eta)=K(\xi,\eta)+\{\Rho(\xi),\eta\}-\{\Rho(\eta),\xi\}.
$$
In particular, $W^{(i)}=K^{(i)}$ for all $i\leq 1$. 
\end{prop}
\begin{proof}
Let $\xi$ be a section of $\Cal A_i$ and $\eta$ be a section of $\Cal
A_j$, with $i,j<0$. To compute $W(\xi,\eta)$, we first have to view
$\xi$ and $\eta$ as vector fields on $M$ via $\tau_-:TM\cong\Cal
A_-$. Then, by construction the section $W(\xi,\eta)$ of $\Cal A$
corresponds to the function $\Cal G_0\to\frak g$ given by
$$
d\tau(\xi^h,\eta^h)+[\tau(\xi^h),\tau(\eta^h)],
$$
where the subscript $h$ denotes the horizontal lift with respect to
the principal connection $\tau_0$. Thus, the $\frak g_0$--components
of $\tau(\xi^h)$ and $\tau(\eta^h)$ are automatically zero, so we may
write 
\begin{align*}
[\tau(\xi^h),\tau(\eta^h)]=\ &[\tau_-(\xi^h),\tau_-(\eta^h)]+
[\tau_+(\xi^h),\tau_-(\eta^h)]+
\\&[\tau_-(\xi^h),\tau_+(\eta^h)]+
[\tau_+(\xi^h),\tau_+(\eta^h)]. 
\end{align*}
On the other hand, from the definition of the curvature it is clear,
that the section $K(\xi,\eta)$ corresponds to the function 
$$
d\tau(\xi^h,\eta^h)+[\tau_-(\xi^h),\tau_-(\eta^h)]+
[\tau_+(\xi^h),\tau_+(\eta^h)].
$$
Now $\tau_+(\xi^h)$ is exactly the function corresponding to
$\Rho(\xi)$, while $\tau_-(\eta^h)$ is the function corresponding to
$\eta$. (Actually, by construction $\tau_-(\eta)$ has values in $\frak
g_j$ only, but this is not important here.) Since the algebraic
bracket $\{\ ,\ \}$ is simply induced by the Lie bracket on $\frak g$,
the formula for $W(\xi,\eta)$ follows immediately. 

To see the second statement, one just has to notice that the algebraic
bracket is by definition homogeneous of degree zero, while all
nonzero homogeneous components of $\Rho$ have degree at least two. 
\end{proof}

\subsection{Remark}\label{4.13}
Note that the latter result, together with the formula \eqref{Weyltransf} 
for the change of the
Weyl--curvature of a Weyl--structure from \ref{4.6} and the formula
\eqref{Rhotransf}
for the change of the Rho--tensor from \ref{3.9}, gives us a formula
for the change of the total curvature of a Weyl--structure under the change
of the Weyl--structure.

\subsection{The construction of normal Weyl--forms}\label{4.14}
Now we are ready to describe the procedure of step by step reducing
the condition of normality of a Weyl--form $\tau\in\Om^1(\Cal
G_0,\frak g)$ to a condition on $\tau_{\leq 0}$ and at the same time
computing step by step the Rho--tensor. From Proposition \ref{4.12} we
know that $W^{(1)}=K^{(1)}$ and from \ref{4.6} we know that this is
actually the same expression for any normal Weyl--form. Usually, this
can be computed in advance, and thus gives us a condition on the
restriction of $\tau_i$ to $T^{i-1}\cg_0$ for $i\leq 0$. Next, by
Proposition \ref{4.12}, we have 
\begin{align*}
W^{(2)}(\xi,\eta)&=K^{(2)}(\xi,\eta)+\{\Rho^{(2)}(\xi),\eta\}-
\{\Rho^{(2)}(\eta),\xi\}\\
&=(K^{(2)}-\partial\Rho^{(2)})(\xi,\eta). 
\end{align*}
If $W^{(2)}$ is to be $\partial^*$--closed, then this implies that
$\partial^*(K^{(2)})=\partial^*\partial\Rho^{(2)}$. On the other hand,
since $H^1(\frak g_-,\frak g)$ is concentrated in homogeneous degrees
less or equal to one and $H^0(\frak g_-,\frak g)=\frak g_{-k}$,
the Hodge decomposition implies that
$\Rho^{(2)}=\square^{-1}\partial^*\partial\Rho-\partial\al_2$ for a
unique smooth section $\al_2$ of $\Cal A_2$. Moreover, since
$\Rho^{(2)}$ has to have values in $\Cal A_+$, it follows that the
restriction of $\square^{-1}\partial^*(K^{(2)})$ to $\Cal
A_{-k}\oplus\dots\oplus\Cal A_{-2}$ must be given by
$\partial(\al_2)$, which gives a condition on the restriction of
$\tau_i$ to $T^{i-2}\Cal G_0$ for $i\leq 0$. If this is satisfied,
then $\al_2$ is uniquely determined, and we can compute $\Rho^{(2)}$
as $\square^{-1}\partial^*(K^{(2)})-\partial\al_2$. Let us notice, how
simple the latter step gets for $|1|$--graded examples: then there is no
$\alpha_2$, the entire forms $\Rho$ and $K_0$ are of homogeneous degree two, 
and so $\Rho$ is simply obtained in the unique way by the formula 
$\Rho=\square^{-1}\partial^*K_0$ promised in Example \ref{4.8}. 

Now this process can be easily iterated. We next consider $K^{(3)}$
which depends only on the (known) component $\Rho^{(2)}$ of the
Rho--tensor and on the restrictions of $\tau_i$ to $T^{i-3}\Cal G_0$
for $i\leq 0$. As above, the restriction of
$\square^{-1}\partial^*(K^{(3)})$ to $\Cal A_{-k}\oplus\dots\oplus\Cal
A_{-3}$ must be given by $\partial(\al_3)$ for a section $\al_3$ of
$\Cal A_3$, which gives conditions on the restrictions of $\tau_i$ for
$i\leq 0$. If these are satisfied, $\al_3$ is uniquely determined, and
we can compute $\Rho^{(3)}$. Finally, once we have reached $K^{(k)}$,
there are no more conditions, since $\tau_{\leq 0}$ is already
completely determined at this stage, so we only get a way to compute
the remaining homogeneous components of the Rho--tensor.

\end{document}